\documentclass[10pt]{article}
\usepackage{amsfonts}
\usepackage{bbm}
\usepackage{mathrsfs}
 \input amssymb.sty
 \usepackage{mathrsfs,amsmath,amssymb}
\usepackage[dvips]{color}
 \usepackage{amssymb}
 \usepackage[dvips]{color}
\allowdisplaybreaks
 \renewcommand{\baselinestretch}{1}
 
 \newtheorem{remark}{\noindent\mbox{Remark}}[section]
 
 \newtheorem{lemma}{\noindent\mbox{Lemma}}[section]
 \newtheorem{theorem}{\noindent\mbox{Theorem}}[section]
 \newtheorem{proposition}{\noindent\mbox{Proposition}}[section]
 \newtheorem{corollary}{\noindent\mbox{Corollary}}[section]
 
\def\bpim{\overline{M}_{0}\overline{M}_{-1}\cdots
\overline{M}_{-n+1}}
\def\pim{M_0M_{-1}\cdots M_{-n+1}}
 
 \def\bq{\begin{equation}}
 \def\eq{\end{equation}}
 \def\eqn{\end{eqnarray}}
 \def\bqn{\begin{eqnarray}}
 \def\proof{\noindent{\it Proof.~~}}
 \def\qed{\hfill$\Box$\medskip}
 \def\rto{\rightarrow\infty}
 \def\z{\left}
 \def\y{\right}
 \def\no{\nonumber}
 \def\mbe{\mathbb{E}}
 \def\mbp{\mathbb{P}}
 \def\ka{{\kappa}}
 
\begin{document}
\noindent{\LARGE\bf Branching structure for an (L-1) random walk in
random environment  and its applications} 

\noindent{\normalsize   Wenming Hong, Huaming Wang

\vspace{0.2 true cm}
\renewcommand{\baselinestretch}{1.3}\baselineskip 12pt
\noindent{\footnotesize\rm  School of Mathematical Sciences \&  Key
Laboratory of Mathematics and Complex Systems, Beijing Normal
University, Beijing 100875, P.R.China;\\%
 Email: wmhong@bnu.edu.cn, huamingking@mail.bnu.edu.cn
\vspace{4mm}}%

\begin{center}
\begin{minipage}[c]{13cm}
\begin{center}\textbf{Abstract}\end{center}
By decomposing the random walk path, we construct a multitype
branching process with immigration in random environment for
corresponding random walk with bounded jumps in random environment.
Then we give two  applications of the branching structure. Firstly,
we specify the explicit invariant density by a method different with
the one used in Br\'{e}mont \cite{br02} and reprove the law of large
numbers of the random walk by a method known as `` the environment
viewed from  particles". Secondly, the branching structure enables
us to prove a stable limit law, generalizing the result of
Kesten-Kozlov-Spitzer \cite{kks75} for the nearest random walk in
random environment. As a byproduct, we also prove that the total
population of a multitype branching process in random environment
with immigration before the first regeneration  belongs to the
domain of attraction of some $\ka$-stable law.
\vspace{0.1cm}\\
\mbox{}\textbf{Keywords:}\quad random walk, branching process, random environment, stable law.\\
\mbox{}\textbf{AMS Subject Classification}:  Primary 60K37;
secondary 60F05.
\end{minipage}
\end{center}
{\center\section{\hspace{-0.5cm.}  Introduction}\label{intr}}
 Random Walk in
Random Environment (RWRE for short) has been extensively studied
(see e.g. Sznitman \cite{S02} or  Zeitouni \cite{ze04}  for a
comprehensive survey), and has wide range of applications both in
probability theory  and physics, for example, in metal physics and
crystallography (see Hughes \cite{hug96} for an introduction). Two
kinds of randomness are involved in RWRE: first the transition
probability chosen randomly at each state position (called random
environment); and second the random walk, a time homogeneous Markov
chain driven by the chosen transition probability.

Random walk  in random environment with bounded jumps ((L-R) RWRE,
that is, the walk, in every unit of time, jumping no more than $L$
to the left and no more than $R$ to the right, where $R$ and $L$ are
positive integers) was first introduced in Key \cite{key84}. Further
developments can be found in Br\'{e}mont \cite{br02,br04} and
L\"{e}tchikov \cite{let89,let94}. We mention here that (L-R) RWRE is
a special case of random walk in random environment on a strip. RWRE
on a strip was first introduced
 in Bolthausen-Goldsheid \cite{bg00}, where the authors provided a
 criteria for the recurrence and transience of the walk. For further
 development of RWRE on a strip, one can refer to
 Goldsheid \cite{gol07,gol08} and Roitershtein \cite{roi08}.

Branching structure played an important role in proving the limit
properties for the nearest neighborhood RWRE. When the walk is
transient to the right, a branching structure was found in
Kesten-Kozlov-Spitzer \cite{kks75} and an elegant stable limit law
was obtained; and the renewal theorem was proved also relying on the
branching structure (Kesten \cite{K77}). Those fine results were
proved  because the steps of the walk could be calculated accurately
under the branching structure. By using the branching structure,
Alili \cite{A99} (see also Zeitouni \cite{ze04}) got the invariant
density  and consequently proved the Law of Large Numbers (LLN for
short) by `` the environment viewed from  particles", a method which
goes back to Kozlov \cite{K85}.

 Br\'{e}mont has also systematically studied (L-1) RWRE in $\cite{br02},$  where the recurrence and transience,
LLN and some central limit theorem have been proved. One of the main
purpose of this paper is to prove a stable limit theorem  to
generalize Kesten-Kozlov-Spitzer \cite{kks75} which dealt the
nearest setting.

To get a stable limit theorem for the (L-1) RWRE,
 one of the crucial steps is to derive the branching structure. In this paper, we will formulate a multitype
branching process for the walk transient to the right. However, it
is much more complicated than the nearest setting, for (1) There are
overlaps between different steps, that is, there may be jumps down
from $i$ to $i-1,i-2,...,i-L.$ Consequently (2) one could not use
jumps down from $i$ directly as the number of $(n-i)$-th generation
of the branching process any more, because  one cannot figure out
the exact parents of particles in $(n-i)$-th generation.

The idea to deal with this difficulty is to imagine that a jump of
size $l$ down from $i$ by the walk  can be remembered by location
$i-1,i-2,...,i-l.$ In this way we can construct a Multitype
Branching Process in Random Environment with one type-1 Immigrant in
each generation (MBPREI for short) to analyze $T_n,$ the first
hitting time of state $n>0$.

After specifying the corresponding MBPREI, we give two applications.
Firstly we can figure out the invariant density directly from the
branching structure  and avoid introducing the \textbf{(IM)}
condition in Br\'{e}mont \cite{br02}. Consequently we can reprove
directly the LLN for the (L-1) RWRE by a method known as `` the
environment viewed from  particles"; secondly we  prove a stable
limit Theorem for the (L-1) RWRE, generalizing Kesten-Kozlov-Spitzer
\cite{kks75} for the nearest one.

 We now define precisely the model of interests to us.
\subsection{Description of the model}
 Let
$\Lambda=\{-L,...,1\}/\{0\}, $ $\Sigma=\{(x_l)_{l\in\Lambda}\in
\mathbb{R}^{L}:\sum_{l\in\Lambda}x_l=1, x_l\ge 0, l\in\Lambda\}$ the
simplex in $ \mathbb{R}^{L+1},$ and $\Omega=\Sigma^\mathbb{Z}.$ Let
$\mu$ be a probability measure on $\Sigma$ and
$\omega_0=(\omega_0(z))_{z\in\Lambda}$ be a $\Sigma$-valued random
vector with distribution $\mu,$ satisfying
$\sum_{z\in\Lambda}\omega_0(z)=1.$  Let
$\mathbb{P}=\otimes^\mathbb{Z}\mu$ on $\Omega$ making $\omega_x,
x\in\mathbb{Z}$ i.i.d. and satisfying
\begin{equation}\label{c2}
  \mbp(\omega_0(z)/\omega_0(1)\ge \varepsilon, \forall
  z\in\Lambda)=1 \text{ for some }\varepsilon>0.
\end{equation}
 Define the shift operator
$\theta$ on $\Omega$ by the relation
\begin{equation}\label{sh}
(\theta\omega)_i=\omega_{i+1}.
\end{equation}
 The pair $(\Omega, \mbp)$ will serve as
the space of environment for both the random walk with bounded jumps
and the multitype branching process which we will give next. The
random walk in random environment $\omega$ with bounded jumps is the
Markov chain defined by $X_0=x$ and transition probabilities
\begin{equation}
    P_{x,\omega}(X_{n+1}=y+z|X_n=y)=\omega_y(z), \forall y\in
    \mathbb{Z}, z\in \Lambda.\no
\end{equation}
In the sequel we refer to $P_{x,\omega}(\cdot)$ as the ``quenched"
law. One also defines the ``annealed" laws
\begin{equation}
  P_x:=\mathbb{P}\times P_{x,\omega}\mbox{\ \ for\ }x\in
  \mathbb{Z}.\no
\end{equation}
  In the rest of the paper, we use respectively
$\mathbb{E}$ corresponding to $\mathbb{P}$, $E_{x,\omega}$
corresponding to $P_{x,\omega}$ and $E_x$ corresponding to $P_x$ to
denote the expectations. And for simplicity, we may use $P,$
$P_\omega,$ $E,$ and $E_\omega$ instead of $P_0,$ $P_{0,\omega},$
$E_0$ and $E_{0,\omega}$.
\subsection{Notations, basic conditions and known results}
All the vectors, both row vectors and column vectors, involved in
this paper are in $\mathbb{R}^{L}$  except otherwise stated. All the
matrices involved are in $\mathbb{R}^{L\times L}.$ Let $x $ be a
vector and $M$ be a matrix. We put
$$|x|=\sum_{i=1}^L|x_i|,\ \|M\|_c=\max_{|y|=1}|My|,\ \text{and } \|M\|=\max_{|x|=1}|xM|.$$
One should note that the matrix norms $\|\cdot\|$ and $\|\cdot\|_c$
are different.

Let $S_{L-1}:=\{x\in \mathbb{R}^L:|x|=1\}$ being the unit ball in
$\mathbb{R}^L,$ and $S_+=\{x\in S_{L-1}:x_i\ge 0\}.$

 For $i\in \mathbb{Z}$ let
\begin{equation}\label{mi}
M_i=\left(
             \begin{array}{cccc}
               b_{i}(1) & \cdots & b_{i}(L-1) & b_{i}(L) \\
               1+b_{i}(1) & \cdots & b_{i}(L-1) & b_{i}(L) \\
               \vdots & \ddots & \vdots & \vdots \\
                b_{i}(1) & \cdots &1+ b_{i}(L-1) & b_{i}(L)\\
             \end{array}\right),\ \overline{M}_i=\left(
                                           \begin{array}{cccc}
                                             a_i(1) & \cdots & a_i{(L-1)} & a_i(L) \\
                                             1 & \cdots & 0 & 0 \\
                                             \vdots & \ddots & \vdots & \vdots \\
                                             0 & \cdots & 1 & 0\\
                                           \end{array}
                                         \right),
\end{equation}
with $a_i(l)=\frac{\omega_i(-l)+\cdots+\omega_i(-L)}{\omega_i(1)},
b_i(l)=\frac{\omega_i(-l)}{\omega_{i}(1)}, 1\le l\le L.$

 We also introduce the following
special vectors. $e_i$ is the unit row vector with $i$-th component
being $1.$ The black $\textbf{1}=(1,...,1),$
$e_0=(\frac{1}{L},...,\frac{1}{L}),$ $x_0=(2,1,...,1)^T,$ and
$\overline{x}_0=(2,-1,0,...,0)^T.$

For $n\ge 0,$ define
$$T_n=\inf[k\ge 0: X_k=n].$$ Note that $T_n$ is the first time the
walk reaching $n.$

Let $\rho=\frac{1-\omega_0(1)}{\omega_0(1)}.$ We use the following
conditions in the paper.

\noindent\textbf{Condition C}\\
    \textbf{(C1)} $\mbe\z(\log^+\rho\y)<\infty,$ with $\log^+x:=0\vee\log x.$\\
     \textbf{(C2)}  $\mathbb{P}\z(\rho>1\y)>0.$ \\
    Note that under \textbf{(C2)} it is an easy task to show that there
    exists some $\kappa_0>2$ such that
    \begin{equation}\label{mkp}
\mathbb{E}\Big[\Big(\min_{1\le i\le L}\{\sum_{j=1}^L
M_0(i,j)\}\Big)^{\kappa_0}\Big]=\mbe(\rho^{\ka_0})> 1.
  \end{equation}
  We mention that (\ref{mkp}) corresponds to (1.13) of Kesten
  \cite{kes73}. But they have different forms since the norm used here differs from the one used in Kesten \cite{kes73}.
  Now fixing such $\kappa_0$ we give a new condition\\
    \textbf{(C3)} $\mathbb{E}(\rho^{\kappa_0}\log^+\rho)<\infty.$\\
Let $\varrho$ be the greatest eigenvalue of $ M_0$. One follows from
(\ref{c2})
that $\varrho>0.$\\
 \textbf{(C4)} The group generated by
$\mbox{supp}[\log\varrho]$ is dense in $\mathbb{R}.$\qed

In the remainder of the paper, except otherwise stated, we always
assume that Condition C holds.
\begin{remark}
The above conditions look more or less like the conditions of the
theorem in Kesten-Kozlov-Spitzer \cite{kks75}. (C1) implies that
$\mbe(\log^+\|M_0\|)<\infty,$ enabling us to calculate the Lyapunov
exponents. Condition (C2) excludes the biased trivial case  and it
also ensures the existing of a number $\ka\in(0,\ka_0]$ such that
$$
  \lim_{n\rto}\{\mbe\z(\parallel
  \pim\parallel^{\kappa}\y)\}^{\frac{1}{n}}=1
$$ (see (\ref{ka}) below).
   One follows from (C3) that $$\mbe\z(\|M_0\|^{\ka_0}\log^+\|M_0\|\y)<\infty, \text{and} \max_{1\le l\le
  L}\mbe\z(\z(\omega_0(-l)/\omega_0(1)\y)^{\ka_0}\y)<\infty,$$
  which will be used many times in this paper.  Also, (C3) implies (C1). Condition (C4) is
  the request of Kesten \cite{kes73} for the proof of the renewal theory of the products of random matrices.
\end{remark}
 The (L-1) RWRE has been studied intensively in Br\'{e}mont \cite{br02},
 where the recurrence and transience criteria, the LLN and some central
 limits theorem have been derived. We state  only the
 recurrence and transience of the model here.
Under condition (C1) one can calculate
   the greatest Ljapounov exponents of
$\{\overline{M}_i\}$ and $\{M_i\}$ under both the norms
$\|\cdot\|_c$ and $\|\cdot\|.$ Indeed in Proposition \ref{lja}
below, we show that $\{M_i\}$ and $\{\overline{M}_i\}$ share the
same greatest Lyapunov exponent. Also as $\{M_i\}$ and
$\{\overline{M}_i\}$ to be considered, it causes no difference to
calculate the Lyapunov exponents under different norms $\|\cdot\|_c$
and $\|\cdot\|.$  So in the remainder of the paper, {\it we use
$\gamma_L$ to denote the greatest Lyapunov exponent of both
$\{M_i\}$ and $\{\overline{M}_i$\}.} The number $ \gamma_L$ provides
the criteria for the transience and recurrence of (L-1) RWRE. We
have

\vspace{0.1cm}
 \noindent\textbf{Theorem A (Br\'{e}mont)} {\it The (L-1) RWRE $\{X_n\}_{n\ge0}$ is $P$-a.s. recurrent,
transient to the right or transient to the left according as
$\gamma_L=0,$ $\gamma_L<0$ or $\gamma_L>0.$}
\subsection{Main results}
 When the walk $\{X_n\}$ is transient to the right ($P$-a.s.), we can
formulate a related MBPREI (with {\it negative time}) to calculate
the steps of the walk. We first define an MBPREI $\{Z_{-n}\}_{n\ge
0},$ with negative time.
 For each integer $k$ we define $Z(k,m)$ to be
the $L$-type branching process in random environment which begins at
time $k.$ That is to say, conditioned on $\omega,$
\begin{eqnarray}\label{imgr}
&&P_\omega(Z(k,m)=\boldsymbol{0})=1,\mbox{ if }m>k,\no\\
&&P_\omega(Z(k,k)=e_1)=1,
\end{eqnarray}
and for $m<k$
\begin{eqnarray}
&&P_\omega(Z(k,m)=(u_1,u_2,...,u_L)\big|Z(k,m+1)=e_1)\no\\
&&\quad\quad=\frac{(u_1+u_2+\cdots+u_L)!}{u_1!u_2!\cdots
u_L!}\omega_{m+1}(-1)^{u_1}\omega_{m+1}(-2)^{u_2}\cdots\omega_{m+1}(-L)^{u_L}\omega_{m+1}(1),\label{zp1}\\
&&P_\omega\z(Z(k,m)=(u_1,...,u_{l-2},u_{l-1}+1,u_l,...,u_L)\big|Z(k,m+1)=e_l\y)\no\\
&&\quad\quad =\frac{(u_1+u_2+\cdots+u_L)!}{u_1!u_2!\cdots
u_L!}\omega_{m+1}(-1)^{u_1}\omega_{m+1}(-2)^{u_2}\cdots\omega_{m+1}(-L)^{u_L}\omega_{m+1}(1),\no\\
&& \quad\quad l=2,3,...,L.\label{zp2}
\end{eqnarray}
In addition to assume that conditioned on $\omega,$ each of the
process $Z(k,*)$ has independent lines of descent, we also assume
that conditioned on $\omega,$ the processes $Z(k,*)$ are
independent.

 Let
\begin{equation}\label{net}
   Z_{-n}=\sum_{k=0}^{n-1}Z(-k,-n),\quad n>0.
\end{equation}
$Z_{-n}$ is the total number of offspring, born at time $-n$ to the
immigrants who arrived between $0$ and $-n+1$, of an MBPREI
beginning at time zero.

Next we consider the path of the (L-1) RWRE $\{X_n\}$ with initial
value $X_0=0.$

 Fix $n>0.$ For $-\infty<i<n,$ $1\le l\le L,$ let
 $$U_{i,l}^n=\#\{0<k<T_n:X_{k-1}>i,X_k=i-l+1\}$$ recording all steps
 by the walk between time interval $(0,T_n)$ from above $i$ to $i-l+1.$
Set
$$U^n_i:=(U_{i,1}^n,U_{i,2}^n,\cdots,U^n_{i,L}).$$ One sees that
$|U^n_i|$ is the total number of steps by the walk reaching or
crossing $i$ downward from above $i.$

Let $I_k=:\{X_m:T_{k}\le m<T_{k+1}\}, k=0,1,...,n-1,$ decomposing
the the random walk path before time $T_n$ into $n$ independent and
non-intersecting pieces. Define for $1\le l\le L,$ $i<k,$
$$U^n_l(k,i)=\#\{T_k\le m<T_{k+1}: X_{m-1}>i,X_m=i-l+1\},$$ counting
the  steps in $I_k$ from above $i$ to $i-l+1.$  Let
$$U^n(k,i)=(U^n_1(k,i),U^n_2(k,i),...,U^n_L(k,i))$$
recording all steps in $I_k$ reaching or crossing $i$ downward from
above $i.$ Set $U^n(k,i)=\textbf{0}$ for $i>k$ and set
$U^n(k,k)=e_1.$ One can see from the definitions of $U^n_i$ and
$U^n(k,i)$ that
\begin{equation*}
U_{i}=\sum_{k=(i+1)\vee0}^{n-1}U^n(k,i).
\end{equation*}

 The relationship between the (L-1) RWRE and the MBPREI $\{Z_{-n}\}_{n\ge0}$ is summarized in
the following theorem which will be proved in Section 2.
\begin{theorem}\label{edis} Suppose that $\gamma_L<0$ (implying that $X_n\rightarrow \infty$ $P$-a.s. by Theorem A).
Then one has that

\noindent(a)
\begin{equation}\label{tn}
T_n=n+\sum_{i=-\infty}^{n-1}|U^n_i|+\sum_{i=-\infty}^{n-1}U^n_{i,1}=n+\sum_{i=-\infty}^{n-1}U^n_ix_0;
\end{equation}

\noindent(b) Each of the processes $U^n(k,*),$ $0\le k\le n-1,$ is
an inhomogeneous multitype branching process beginning at time $k$
with branching mechanism
\begin{eqnarray}\label{up01}
&&P_\omega(U^n(k,i-1)=(u_1,...,u_L)\big|U^n(k,i)=e_1)\no\\
&&\hspace{2cm}=\frac{(u_1+\cdots+u_L)!}{u_1!\cdots
u_L!}\omega_i(-1)^{u_1}\cdots\omega_i(-L)^{u_L}\omega_i(0),
\end{eqnarray}
and for $2\le l\le L,$
\begin{eqnarray}\label{up02}
&&P_\omega\z(U^n(k,i-1)=(u_1,...,1+u_{l-1},...,u_L)\big|U^n(k,i)=e_l\y)\no\\
 &&\hspace{2cm}=\frac{(u_1+\cdots+u_L)!}{u_1!\cdots
u_L!}\omega_i(-1)^{u_1}\cdots\omega_i(-L)^{u_L}\omega_i(0).
\end{eqnarray}
Moreover conditioned on $\omega,$ $U^n(k,*),k=0,1,...,n-1$ are
mutually independent and each of the branching processes $U^n(k,*)$
has independent line of descent. Consequently
$U^n_{n-1}=0,U^n_{n-2},\cdots, U^n_{1},U_{0}^n$ are the first $n$
generations of an inhomogeneous multitype branching process with a
type-1 immigration in each generation in random environment.

\noindent(c) $U^n_{n-1}=0,U^n_{n-2},\cdots, U^n_{1},U_{0}^n$ has the
same distribution with the first $n$ generations of the
inhomogeneous MBPREI $\{Z_{-n}\}_{n\le 0}$ defined in (\ref{net}).
\end{theorem}
We have immediately the following corollary about the offspring
matrices of the multitype branching process $\{U^n_i\}_{0\le i\le
n-1}.$
\begin{corollary}\label{col}
For the process $\{U_i^n\}_{i=0}^{n-1},$ let $ M_i$ be the $L\times
L$ matrix whose $l$-th row is the expected number of offspring born
to a type-$l$ parent of the $(n-i)$-th generation conditioned on
$\omega,$ that is, $E_\omega\z(U^n(i,i-1)|U^n(i,i)=e_l\y).$ Then one
has that
\begin{equation*}
 M_i=\left(
             \begin{array}{cccc}
               b_{i}(1) & \cdots & b_{i}(L-1) & b_{i}(L) \\
               1+b_{i}(1) & \cdots & b_{i}(L-1) & b_{i}(L) \\
               \vdots & \ddots & \vdots & \vdots \\
                b_{i}(1) & \cdots &1+ b_{i}(L-1) & b_{i}(L)\\
             \end{array}
           \right)
\end{equation*} which coincides with the definition of $M_{i}$ in
(\ref{mi}). Similarly for the process $\{Z_{-n}\}_{n\ge 0}$ let $
M_{-i}$ be the $L\times L$ matrix whose $l$-th row is the expected
number of offspring born to a type-$l$ parent at time $-i,$
conditioned on $\omega:$ $E_\omega\z(Z(-i,-i-1)|Z(-i,-i)=e_l\y).$
Then one has that
\begin{equation*}\label{mb1}
 M_{-i}=\left(
             \begin{array}{cccc}
               b_{-i}(1) & \cdots & b_{-i}(L-1) & b_{-i}(L) \\
               1+b_{-i}(1) & \cdots & b_{-i}(L-1) & b_{-i}(L) \\
               \vdots & \ddots & \vdots & \vdots \\
                b_{-i}(1) & \cdots &1+ b_{-i}(L-1) & b_{-i}(L)\\
             \end{array}
           \right).
\end{equation*}
\end{corollary}
Part (c) of Theorem \ref{edis} says that
$U^n_{n-1}=0,U^n_{n-2},\cdots, U^n_{1},U_{0}^n$ has the same
distribution with the first $n$ generations of the inhomogeneous
MBPREI $\{Z_{-n}\}_{n\le 0}$ defined in (\ref{net}). Instead of
studying the limiting behaviors of the hitting time $T_n$ directly,
one turns to study that of the $L$-type branching process
$\{Z_{-n}\}_{n\ge 0}$ first.

Let $\nu_0\equiv0,$ and define recursively
\begin{equation}
    \nu_n=\min\{m>\nu_{n-1}:Z_{-m}=\boldsymbol{0}\}\mbox{ for } n>0,\no
\end{equation} being the successive regeneration times of MBPREI $\{Z_{-n}\}_{n\ge0}.$
 For simplicity we write
$\nu_1$ as $\nu.$

One sees that the regeneration time $\nu$ of the MBPREI
$\{Z_{-n}\}_{n\ge0}$ corresponds to the regeneration position (some
position where the walk will never go back after passing it) for
$\{X_n\}_{n\ge0}.$

Define also
\begin{equation}
    W=\sum_{k=0}^{\nu-1}Z_{-k},\no
\end{equation}
the total number of offspring born before the regeneration time
$\nu.$

Kesten\cite{kes73} (see Theorem \ref{kes} below) has proved that if
Condition C  holds and  the greatest Lyapunov exponent $\gamma_L$ of
$ \{M_{-n}\}_{n\ge0}$ is strictly negative,
  there exists a unique $\kappa\in(0,\kappa_0],$ such that
\begin{equation}\label{ka}
  \log\rho(\kappa)=\lim_{n\rto}\frac{1}{n}\log\mbe\z(\parallel \pim\parallel^{\kappa}\y)
  =0.
\end{equation}
Then we have the following limiting theorem of MBPREI
$\{Z_{-n}\}_{\ge0}.$

\begin{theorem}\label{siz0} Let $\ka$ be the number in (\ref{ka}).
Suppose that   Condition C holds and $\gamma_L<0.$  If  $\kappa> 2,$
then $ E((Wx_0)^2)<\infty;$ if $\kappa\le 2,$ then there exists some
$0<K_3<\infty$ such that
\begin{equation}\label{ws0}
\lim_{t\rto}t^\kappa P(Wx_0\ge t)=K_3.
\end{equation}
\end{theorem}

  For $ n\ge 0$ define
$\overline{\omega}(n)=\theta^{X_n}\omega.$ Then
$\{\overline{\omega}(n)\}$ is a Markov chain with transition kernel
$$\overline{P}(\omega,d\omega')=\omega_0(1)\delta_{\theta\omega=\omega'}+\sum_{l=1}^L\omega_0(-l)\delta_{\theta^{-l}\omega=\omega'}.$$ In \cite{br02} an \textbf{(IM)} condition is
said to be satisfied if there is  $\pi(\omega)$ such that
$$\int \tilde{\pi}(\omega)\mbp(d\omega)=1\mbox{ and }
\tilde{\pi}(\omega)=\overline{P}*\tilde{\pi}(\omega),$$ where
$\tilde\pi(\omega)=\pi(\omega)[\mathbb{E}(\pi(\omega))]^{-1}.$ Under
\textbf{(IM)} condition Br\'{e}mont showed an LLN of $\{X_n\}$ in
\cite{br02}. But the \textbf{(IM)} condition was not given directly
in the words of environment $\omega.$ So one has to check the
existence of the invariant density $\pi(\omega).$ In \cite{br02},
Br\'{e}mont showed the existence of $\pi(\omega)$ by analyzing its
definition and the transition probability of the walk.

What makes difference in this article is that, with the help of the
branching structure, we specify the invariant density $\pi(\omega)$
directly by analyzing a multitype branching process. Therefore we
avoid introducing the \textbf{(IM)} condition and show directly that
$\{X_n\}$ satisfies an LLN with a positive speed under the
assumption $``\mathbb{E}(\pi(\omega))<\infty"$. Also the speed has
 a simple explicit form $[\mathbb{E}(\pi(\omega))]^{-1}.$

Define $\pi(\omega):=\frac{1}{\omega_0(1)}\z(1+\sum_{i=1}^\infty
e_1\overline{M}_i\cdots \overline{M}_1e_1^T\y).$ Let
$\tilde{\pi}(\omega)=\frac{\pi(\omega)}{\mbe(\pi(\omega))}.$ Then we
have
\begin{theorem}\label{lln}
  Suppose that $\mbe(\pi(\omega))<\infty.$
  Then we have that\\
  \mbox{\quad\quad} (i)  $\gamma_L<0;$\\
  \mbox{\quad\quad} (ii) $\tilde{\pi}(\omega)\mathbb{P}(d\omega)$ is invariant under the
  kernel $\overline{P}(\omega, d\omega'),$ that is
\begin{equation}\label{invd}
\int1_B\tilde{\pi}(\omega)\mathbb{P}(d\omega)=\iint 1_{\omega'\in
B}\overline{P}(\omega,d\omega')\tilde{\pi}(\omega)\mathbb{P}(d\omega);
\end{equation}
\mbox{\quad\quad} (iii) and
  $\mbp$-a.s.,
  $\lim_{n\rto}\frac{X_n}{n}=\frac{1}{\mbe(\pi(\omega))}.$
\end{theorem}
\begin{remark}
  The
  independence assumption of the environment is unnecessary. It is
  enough if $(\Omega, \mbp, \theta)$ is an ergodic system.  \end{remark}
  Also for the (L-1) RWRE we have the following stable limit
  theorem, generalizing  Kesten-Kozlov-Spitzer
\cite{kks75} which dealt with the nearest setting.
\begin{theorem}\label{main}
Suppose that Condition C holds and that $\gamma_L<0.$ Let $\ka$ be
the number in (\ref{ka}).  Let $L_\ka(x)$ be a $\ka$-stable law
($L_\ka$ is concentrated on $[0,\infty)$ if $\ka<1$ and has mean
zero if $\ka>1$). Then with $0<A_\ka, B_i<\infty$ suitable
constants, $\Phi(x):=\frac{1}{\sqrt{2\pi}}\displaystyle
\int_{-\infty}^x
e^{-\frac{s^2}{2}}ds,$\\
 (i) if $0<\ka<1,$ \begin{align}
&\lim_{n\rto}P(n^{-\frac{1}{\ka}}T_n\le
  x)=L_\ka(x),\no\\
  &\lim_{n\rto}P(n^{-\ka}X_n\le
  x)=1-L_\ka(x^{-\frac{1}{\ka}});\no
\end{align}
  (ii) if $\ka=1,$ then for suitable $D(n)\sim \log n$ and $\delta(n)\sim (A_1\log n)^{-1}n,$
  \begin{align*}
&\lim_{n\rto}P(n^{-1}(T_n-A_1nD(n\mu^{-1}))\le x)=L_1(x),\\
  &\lim_{n\rto}P(n^{-1}(\log n)^2(X_n-\delta(n))\le
  x)=1-L_1(-A_1^2x);
  \end{align*}
  (iii) if $1<\ka<2,$
  \begin{align*}
&\lim_{n\rto}P\z(n^{-\frac{1}{\ka}}(T_n-A_\ka n)\le
  x\y)=L_\ka(x),\\
& \lim_{n\rto}P\z(n^{-\frac{1}{\ka}}\z(X_n-\frac{n}{A_\ka}\y)\le
  x\y)=1-L_\ka(-xA_\ka^{1+\ka^{-1}});
  \end{align*}
(iv) if $\ka=2,$
\begin{align*}
&\lim_{n\rto}P\Big(\frac{T_n-A_2 n}{B_1\sqrt{n\log n}}\le
x\Big)=\Phi(x),\\
&\lim_{n\rto}P\z(A_2^{\frac{3}{2}}B_1^{-1}(n \log
n)^{-\frac{1}{2}}\z(X_n-\frac{n}{A_2}\y)\le x\y)=\Phi(x);
\end{align*}
(v) if $\ka>2,$ \begin{align*} &\lim_{n\rto}P\Big(\frac{T_n-B_3
n}{B_2\sqrt{n}}\le x\Big)=\Phi(x),\\
&\lim_{n\rto}P\z(B_3^{\frac{3}{2}}B_2^{-1}n^{-\frac{1}{2}}\z(X_n-\frac{n}{B_3}\y)\le
x\y)=\Phi(x).
\end{align*}
\end{theorem}
\textbf{Notes:} The stable limit law for the nearest neighborhood
RWRE ((1-1) RWRE) was shown in Kesten-Kozlov-Spitzer \cite{kks75}.
But to prove Theorem \ref{main} is far more than a trivial work for
the following reasons:
\begin{itemize}
  {\rm\item[(1)]The branching structure (MBPREI $\{Z_{-n}\}_{n\ge0}$) for (L-1) RWRE was never seen in literatures we are aware of. But it is crucial to construct such branching structure to prove Theorem
  \ref{main}.
  \item[(2)] After constructing the branching structure, to prove Theorem \ref{main}, a key step is to
  show Theorem \ref{siz0}, that is, to show that $W,$ the total population of $\{Z_{-n}\}_{n\ge0}$ before the first regeneration,   belongs to the domain of attraction of a $\ka$-stable
  law:
    \begin{equation}\label{ws1}
\lim_{t\rto}t^\kappa P(Wx_0\ge t)=K_3.
  \end{equation} For this purpose one should use the
  tail of the series of the products of random matrices of Kesten \cite{kes73}, that is,
  \begin{equation}\label{fin}
\lim_{t\rto}t^\ka \mathbb{P}\z(\sum_{n=0}^\infty x\pim
x_0\y)=K(x,x_0),
  \end{equation} where $K$ is a constant depending on  positive $x\in \mathbb{R}^L.$

But to prove (\ref{ws1}) one needs to find out how the constant $K$
depending on $x$ explicitly. For general random matrix, this is
still open. But for $\{M_i\}$, by the similarity between $M_i$ and
$\overline{M}_i,$ we prove a finer result based on (\ref{fin}). We
show that
\begin{equation*} \lim_{t\rto}t^\ka \mathbb{P}\z(\sum_{n=0}^\infty
x\pim x_0\y)=K_2|xB|^{\kappa},
  \end{equation*}}
  where $K_2$ is independent of $x.$
\end{itemize}
   For the reason we list above, although Key, E.S. indicated  that
  ``the general argument for finding limiting distributions for $\{X_n\}$ seems to go through except now $\{Z_t\}$ (some branching process) is
  multitype; the only part that seems not to be line by line rewriting of Kesten, Kozlov and Spitzer's argument is the proof of
  \begin{equation}\label{nue}
    P(\nu>t)<K_4\exp(-K_5t)"
  \end{equation}
  $\nu$ being the regeneration time of some multitype branching process $\{Z_t\}$(see \cite{key87} page 350), we think it makes sense to prove Theorem
  \ref{main}.\qed

Since the proof of Theorem \ref{main} will be a long march. We
describe the skeleton of its proof.

 In order to determine the limit law of $X_n,$ we
consider first the limit law of the hitting time $T_n.$  One gets
from Theorem \ref{edis}  that
\begin{equation*}\label{tnul}
    T_n=n+\sum_{i=-\infty}^{n-1}|U^n_i|+\sum_{i=-\infty}^{n-1}U^n_{i,1}=n+\sum_{i=-\infty}^{n-1}U^n_ix_0.
\end{equation*}
 Note that when the walk transient to the right, there are only finite steps in $(-\infty,0),$
 that is, $P$-a.s., $\frac{1}{n}\sum_{i<0}U_i^n x_0\rightarrow 0. $
 Consequently,  it suffices to show that $\sum_{i=0}^{n-1}U_i^n x_0$ converges
 to $L_\ka$ in distribution after suitable normalization.
 Also in Theorem \ref{edis}
 we observe that $$U^n_{n-1}=0,U^n_{n-2},\cdots, U^n_{1},U_{0}^n$$
has the same distribution with the first $n$ generations of an
inhomogeneous MBPREI $\{Z_{-n}\}_{n\ge0}$ such that conditioned on
$\omega,$ $E_\omega(Z_{-t}|Z_{-k},0\le k<t)=(Z_{-t+1}+e_1)
M_{-t+1}.$

Let $W_k:=\sum_{\nu_k\le t<\nu_{k+1}}Z_{-t}$ be the total offspring
born between time interval $[\nu_k,\nu_{k+1}).$ Then due to the
independence of the environment, $(\nu_{k+1}-\nu_k, W_k),\
k=0,1,2,...$ are independent and identically distributed.

The key step is to show that
\begin{equation}\label{dom}
W_0x_0\mbox{ is in the domain of attraction of a $\ka$-stable law,}
\end{equation}  which is proved in Theorem \ref{siz0}. Then Theorem \ref{main}
follows by a standard argument. To show (\ref{dom}), we could
approximate $Wx_0$ by some random variable of the form
$\Gamma(R+I)x_0,$ where $\Gamma$ is a random row vector with
positive components and independent of $R,$ and $Rx_0$ has the same
distribution with
$$\eta_0x_0:=\sum_{m=1}^\infty  M_0M_{-1}\cdots  M_{-m+1}x_0.$$

It remains to show that for all row vector $x$ with positive
components, $x\eta_0x_0$ belongs the domain of attraction of
$\ka$-stable law, that is
\begin{equation}\label{apd}P(x\eta_0x_0>t)\sim K_2 |xB|^\ka
t^{-\ka} \mbox{ as } n\rto
\end{equation}
with the constant $K_2$ independent of $x.$

To this end we first analyze the connections between the matrices
$M_i$ and $\overline{M}_i$ and find that  the projection of
$\eta_{0}x_0$ on different directions $e_l,1\le l\le L,$ that is
$e_l\eta_{0}x_0,$ have the same distributions up to certain linear
transformations (see Proposition \ref{tran} below).

This fine property of the random variable $\eta_0$ together with
Kesten's results of the products of random matrices enables us to
show that the constant $K_2$ of (\ref{apd}) is independent of $x$
(see Theorem \ref{ep} below).

We arrange  the remainder of this paper as follows. In Section
\ref{dec}, we formulate the related branching structure MBPREI and
express the hitting time  $T_n$ by the MBPREI, i.e., Theorem
\ref{edis} is proved. In Section \ref{lnt} we study the connections
between the matrices $M_i$ and $\overline{M}_i$ which will be
important to show (\ref{apd}). In Section \ref{alln}, we give the
 invariant density from the point of branching structure and provide an alternative proof of LLN (Theorem \ref{lln}).
 Finally the long Section \ref{sta} is devoted to studying the tail of $Wx_0$ (Theorem \ref{siz0}) and to proving the
 stable limit law (Theorem \ref{main}).

{\center {\section{ Branching structure and   hitting
times}\label{dec}} } In this section we assume the walk transient to
the right, i.e., $X_n\rightarrow \infty$ $P$-a.s., and always use
notation $J_{s,h}$ to denote a jump (a piece of random walk path)
taken by $\{X_n\}$ from $s$ to $h.$ That is
$$J_{s,h}=\{(X_n,X_{n+1}):X_n=s,X_{n+1}=h\},n\in\mathbb{N}.$$
 We are going to find a multitype branching process (i.e., MBPREI) to analyze $T_n$.
 To see this, for some $1\le l\le L,$ suppose we have a
jump, say $J_{i,i-l},$ by the walk from $i$ downward to $i-l$ before
time $T_n.$ Though throughout this jump, the walk will not stop at
$i-1,i-2,...,i-l+1,$ we can imagine that  it will cross
$i-1,i-2,\cdots,i-l+1,$ and reach $i-l$ at last. For $-\infty<i<n,$
we want to record how much times the walk will cross or reach $i$
downward from above $i$. So for $-\infty<i<n,$ $1\le l\le L,$ we
define
\begin{equation}\label{uil}
    U_{i,l}^n=\#\{0<k<T_n:X_{k-1}>i,X_k=i-l+1\},\no
\end{equation}
being the records of steps by the walk from above $i$ to $i-l+1,$
and let
\begin{equation}
    U^n_i=(U_{i,1}^n,U_{i,2}^n,\cdots,U^n_{i,L}).\no
\end{equation}
Then $|U_i^n|$ is the total number of times the walk crossing or
reaching $i$ downward from above $i$ before $T_n$. In particular
$U_{i,1}^n$ is the total number of jumps taken by the walk downward
from above $i$ which reach $i.$ But every jump taken by the walk
downward must reach some $i.$ From this point of view, the total
number of steps taken by the walk downward before time $T_n$ is
$\sum_{i=-\infty}^{n-1}U_{i,1}^n.$

 On the other hand, for the walk
transient to the right, since we are considering the  (L-1)  RWRE,
every record of the walk reaching or crossing $i$ downward from
above $i$, an individual of $U_i^n$, corresponds to a jump taken by
the walk upward from $i$ to $i+1.$ Therefore the total number of the
jumps taken by the walk upward before time $T_n$ is
$\sum_{i=-\infty}^{n-1}|U_i^n|.$ Then we conclude from the above
discussion that
\begin{equation}
    T_n=n+\sum_{i=-\infty}^{n-1}|U^n_i|+\sum_{i=-\infty}^{n-1}U^n_{i,1}=n+\sum_{i=-\infty}^{n-1}U^n_ix_0,\no
\end{equation}
where $x_0=(2,1,...,1)^T\in \mathbb{R}^L.$ So instead of studying
$T_n$ directly we consider $\{U_i^n\}_{i<n}.$

We first divide the path between $0$ and $T_n$ into $n$ pieces which
do not have intersection.  For $k=0,1,..,n-1,$ define
$$I_k=:\{X_m:T_k\le m<T_{k+1}\},\tau_{k}=\{t:T_k\le t<T_{k+1}\}$$
Then one follows from the strong Markov property that
$$(I_0,\tau_0),(I_1,\tau_1),...,(I_{n-1},\tau_{n-1})$$ are mutually independent
under the quenched probability $P_\omega^0$. We will see next that
each of the pieces $(I_k,\tau_k),$ $0\le k\le n-1,$ corresponds to
an immigration structure.

Now we fix $0\le k\le n-1.$ We want to construct an $L$-type
branching process from the random walk path $I_k.$

Define for $1\le l\le L,$ $i<k,$
$$U^n_l(k,i)=\#\{T_k\le m<T_{k+1}: X_{m-1}>i,X_m=i-l+1\},$$ counting
the  steps in $I_k$ from above $i$ to $i-l+1.$  Let
$$U^n(k,i)=(U^n_1(k,i),U^n_2(k,i),...,U^n_L(k,i))$$
recording all steps in $I_k$ reaching or crossing $i$ downward from
above $i.$

Note that if $i\ge k,$ there is no step between time interval
$\tau_k$ by the walk reaching or crossing $i$ downward from above
$i$. So we set very naturally $U^n(k,i)=0$ for $i>k.$

But note also that  there may be some  steps between time interval
$\tau_k$ by the walk reaching or crossing $k-1$ down ward from
above. If we want to consider these jumps as the particles of a
branching process at time $k-1,$ we must figure out their parents.
So we can assume that there is an immigrant entering the system at
time $k.$ Therefore we set  $U^n(k,k):=e_1,$ representing the
immigrant entering at time $k.$

We show next that $\{U^n(k,i)\}_{i\le k}$ is an inhomogeneous
$L$-type branching process beginning at time $k.$

 Fix $i<k.$ Let $\eta^k_{i,0}=T_k.$ Define recursively for $j\ge 1$
$$\eta_{i,j}^k=\min\{\eta_{i,j-1}^k<m<T_{k+1}:X_{m-1}>i,X_m\le
i\}.$$
  Then by definition $\eta^k_{i,j},j\ge 1$ are the
successive time of steps in $I_k$ taken by  the walk crossing or
reaching $i$ downward from above $i$. For $i\le k-1$, $1\le l\le L,$
we have
\begin{equation}\label{rela}
U_l^n(k,i)=\sum_{j=1}^\infty
1_{[T_k<\eta^k_{i,j}<T_{k+1},X_{\eta^k_{i,j}}=i-l+1]}.
\end{equation} Also define
$$\xi^{k,j}_{i,l}=\#\{\eta^k_{i+1,j}\le m<\eta^k_{i+1,j+1}<T_{k+1}:
X_{m-1}>i,X_{m}= i-l+1\},$$ recording the steps by the walk from
above $i$ to $i-l+1$ between the $j$-th and the $j+1$-th excursions
reaching or crossing $i+1$ in the time interval $\tau_k,$ and define
$$\xi_i^{k,j}=(\xi_{i,1}^{k,j},...,\xi_{i,L}^{k,j}).$$

 Then it follows from the path decomposition
of $I_k$ and the strong Markov property that $\xi_{i}^{k,j},\
j=1,2,...$ are i.i.d. for fixed $k$ and $i.$

In the definition of $\xi_{i,l}^{k,j},$ where things get delicate is
the first step in the time interval
$[\eta^k_{i+1,j},\eta^k_{i+1,j+1}).$ Note that on the event
$\{X_{\eta^k_{i+1,j}}=i+1\}$ there is no jump of the form $J_{s,h}$
with $s>i+1$  and $h\le i$ during the time interval
$[\eta^k_{i+1,j},\eta^k_{i+1,j+1}).$ Hence for $1\le l\le L$
\begin{eqnarray}
\xi^{k,j}_{i,l}\hspace{-0.2cm}&=&\hspace{-0.2cm}\#\{\eta^k_{i+1,j}\le
m<\eta^k_{i+1,j+1}<T_n:
X_{m-1}>i,X_{m}= i-l+1\}\no\\
\hspace{-0.2cm}&=&\hspace{-0.2cm}\#\{\eta^k_{i+1,j}<
m<\eta^k_{i+1,j+1}<T_n:
X_{m-1}=i+1,X_{m}= i-l+1\}\no\\
\hspace{-0.2cm}&=&\hspace{-0.2cm}\#\{ \mbox{jumps of the kind }
J_{i+1,i+1-l} \mbox{ by the walk during time interval
}(\eta^k_{i+1,j},\eta^k_{i+1,j+1})\}.\no
\end{eqnarray} For $2\le m\le L,$ on the event
$\{X_{\eta^k_{i+1,j}}=i-(m-1)+1\},$ the first step in the time
interval $[\eta^k_{i+1,j},\eta^k_{i+1,j+1})$ taken by the walk from
the above of $i+1$ reaches or crosses $i$ downward and hence
contributes unconditionally one particle to $\xi_{i,m-1}^{k,j}.$
There is  no other step of the form $J_{s,h}$ with $s>i+1$ and $h\le
i$  in such interval. Hence
$$\xi^{k,j}_{i,m-1}=1+\#\{
\mbox{jumps of the kind } J_{i+1,i-(m-1)+1}\mbox{ by the walk during
time interval }(\eta^k_{i+1,j},\eta^k_{i+1,j+1})\},$$ and for $2\le
l\le L, l\neq m-1$
$$\xi^{n,j}_{i,l}=\#\{
\mbox{jumps of the kind } J_{i+1,i+1-l}\mbox{ by the walk during
time interval }(\eta^k_{i+1,j},\eta^k_{i+1,j+1})\}.$$
  Then it follows from the above discussion and the definition of
  $\{U^n(k,i)\}_{i\le k}$ that
\begin{eqnarray}\label{up1}
&&P_\omega(U^n(k,i-1)=(u_1,...,u_L)\big|U^n(k,i)=e_1)\no\\
&&\hspace{1cm}=P_\omega(\xi_{i-1,m}^{k,j}=u_m,1\le m\le
L\big|X_{\eta^k_{i,j}}=i)\no\\
 &&\hspace{1cm}=\frac{(u_1+\cdots+u_L)!}{u_1!\cdots
u_L!}\omega_i(-1)^{u_1}\cdots\omega_i(-L)^{u_L}\omega_i(0),
\end{eqnarray}
and for $2\le l\le L,$
\begin{eqnarray}\label{up2}
&&P_\omega(U^n(k,i-1)=(u_1,...,1+u_{l-1},...,u_L)\big|U^n(k,i)=e_l)\no\\
&&\hspace{1cm}=P_\omega(\xi_{i-1,l-1}^{k,j}=u_{l-1}+1,\xi_{i-1,m}^{k,j}=u_m, m\neq l-1, 1\le m\le L\big|X_{\eta^k_{i,j}}=i+1-l)\no\\
 &&\hspace{1cm}=\frac{(u_1+\cdots+u_L)!}{u_1!\cdots
u_L!}\omega_i(-1)^{u_1}\cdots\omega_i(-L)^{u_L}\omega_i(0).
\end{eqnarray}
Then one can conclude that $\{U^n(k,i)\}_{i\le k}$ is an
inhomogeneous $L$-type branching process beginning at time $k.$

It follows from the independence of the
$\{(I_k,\tau_k)\}_{k=0}^{n-1}$ that $U^n(k,*),k=n-1,...,1,0$ is
mutually independent. Also by the above discussion of path
decomposition inner the random walk piece $I_k,$ we found that
particles of $U^n(k,i)$ generate offspring  independently.

The above two kinds of independence correspond to the independence
we imposed on the branching processes $Z(k,*).$

 Note that by the definition of $U^n_{i,l},$ for  $1\le l\le L,$ one has that
\begin{equation}\label{relb}
U_{i,l}^n=\sum_{k=(i+1)\vee 0}^{n-1}\sum_{j=1}^\infty
1_{[\eta^k_{i,j}<T_n,X_{\eta^k_{i,j}}=i-l+1]}.
\end{equation}
Then it follows from (\ref{rela}) and (\ref{relb}) that
$$U^n_{i,l}=\sum_{k=(i+1)\vee 0}^{n-1}U^n_l(k,i),$$ which implies that
\begin{equation}\label{brui}
U^n_{i}=\sum_{k=(i+1)\vee 0}^{n-1}U^n(k,i)
\end{equation}
 Comparing
the above (\ref{up1}), (\ref{up2}) with (\ref{zp1}), (\ref{zp2}) in
the definition of $\{Z_{-n}\}_{n\ge0},$ since
$\omega_{n-1},\omega_{n-2}...,\omega_1$ have the same joint
distribution as $\omega_0,\omega_{-1},...,\omega_{-n+1},$ it follows
that
$$U^n_{n-1}=0,U^n_{n-2},..., U^n_{1},U_{0}^n$$
has the same distribution with the first $n$ generations of the
inhomogeneous MBPREI $\{Z_{-n}\}_{n\ge0}$ defined in Section
\ref{intr}.

Summarizing  the discussion above, we obtain the proof of Theorem
\ref{edis}.

 {\center\section{Connections between  matrices $M_{i}$ and $\overline{M}_i$} \label{lnt}}
Recall that in \cite{br02}, the author used the greatest Lyapunov
exponent $\gamma_L$ of $\{\overline{M}_i\}$ to characterize the
recurrence and transience of
 (L-1) RWRE $\{X_n\}$ (see also Theorem A above). One may be curious that in Corollary \ref{col} the expected offspring
 matrices of the branching processes $\{U_i^n\}_{0\le i\le n-1}$ and $\{Z_{-n}\}_{n\ge 0}$
 are matrices in the sequence $\{M_i\}_{i\in Z}.$ From this point of
 view, one may guess that there must be some intrinsic connections
 between the matrices $M_i$ and $\overline{M}_i.$ So the main task of
 this short section is to find some specific relations between
 $\{M_i\}$ and $\{\overline{M}_i\}.$

Firstly, one easily sees that $M_i$ and $\overline{M}_i$ are similar
to each other. Indeed, introduce the deterministic matrix
\begin{equation}\label{bi}  B=\left(
             \begin{array}{cccc}
               1 & &  & \\
               1 & 1 &  & \\
               \vdots & \vdots& \ddots &  \\
                1& 1&\cdots  & 1\\
             \end{array}
           \right), \mbox{ with inverse }
 B^{-1}=\left(
             \begin{array}{cccc}
               1 & &  & \\
               -1 & 1 &  & \\
               & \ddots & \ddots &  \\
                &  &-1 & 1
             \end{array}
           \right)
\end{equation}
the entries in the blank being all zero. Then one has that
\begin{equation}\label{simil}
  \overline{M}_i=B^{-1}M_iB.
\end{equation}
Then the similarity between $M_i$ and $\overline{M}_i$ follows.

One notes that since $\mbp$ makes $\{\omega_i\}_{i\in \mathbb{Z}}$
an i.i.d. sequence, it also makes $\{\overline{M}_i\}_{i\in
\mathbb{Z}}$ and $\{ M_i\}_{i\in \mathbb{Z}}$ two i.i.d. sequences
as well. These two random sequences of matrices are of great
importance to us.  Under condition \textbf{(C1)}  one can apply
Oseledec's multiplicative ergodic theorem (see \cite{ose68}) to both
$\{\overline{M}_i\}$ and $ \{M_i\},$ with the shift operator defined
in (\ref{sh}). Write $\gamma_L(\overline{M},\theta)\ge
\gamma_{L-1}(\overline{M},\theta)\ge\cdots\ge
\gamma_1(\overline{M},\theta)$ for the Lyapunov exponents of
$\{\overline{M}_i\}$ and $\gamma_L( M,\theta)\ge \gamma_{L-1}(
M,\theta)\ge\cdots\ge \gamma_1( M,\theta)$ for the Lyapunov
exponents of $ \{M_i\}$ under the matrix norm $\|\cdot\|_c.$ For
simplicity we write $\gamma_L(\overline{M},\theta)$ as
$\gamma_L(\overline{M})$ and $\gamma_L( M,\theta)$ as $\gamma_L(M)$
respectively. Due to the positivity of both $ M_0$ and
$\overline{M}_0,$ we have for all $x\in S_+,$ $\mbp$-a.s.,
\begin{eqnarray}\label{gm}
\gamma_L(\overline{M})&=&\lim_{n\rto}\frac{1}{n}\log\|\overline{M}_{n-1}\cdots
\overline{M}_1\overline{M}_0\|_c=\lim_{n\rto}\frac{1}{n}\log|\overline{M}_{n-1}\cdots \overline{M}_1\overline{M}_0x|\no\\
&=&\lim_{n\rto}\frac{1}{n}\mbe(\log\|\overline{M}_{n-1}\cdots
\overline{M}_1\overline{M}_0\|_c)
\end{eqnarray}
and
\begin{eqnarray}
\gamma_L( M)&=&\lim_{n\rto}\frac{1}{n}\log\| M_{n-1}\cdots
 M_1M_0\|_c=\lim_{n\rto}\frac{1}{n}\log| M_{n-1}\cdots
 M_1M_0x|\no\\
&=&\lim_{n\rto}\frac{1}{n}\mbe(\log\| M_{n-1}\cdots M_1M_0\|_c).\no
\end{eqnarray}
Similarly we can calculate the greatest Lyapunov exponent
$\overline{\gamma}_L(M)$ of $\{M_i\}_{i\le 0}$ and
$\overline{\gamma}_L(\overline{M})$ of $\{\overline{M}_i\}_{i\le0}$
under matrix norm $\|\cdot\|$ as
\begin{eqnarray*}
\overline\gamma_L(\overline{M})&=&\lim_{n\rto}\frac{1}{n}\log\|\overline{M}_{0}\overline
M_{-1}\cdots
\overline{M}_{-n+1}\|=\lim_{n\rto}\frac{1}{n}\log|x\overline{M}_{0}\overline
M_{-1}\cdots
\overline{M}_{-n+1}|\no\\
&=&\lim_{n\rto}\frac{1}{n}\mbe(\log\|\overline{M}_{0}\overline
M_{-1}\cdots \overline{M}_{-n+1}\|)
\end{eqnarray*}
and
\begin{eqnarray*}
\overline\gamma_L( M)&=&\lim_{n\rto}\frac{1}{n}\log\| M_{0}
M_{-1}\cdots M_{-n+1}\|=\lim_{n\rto}\frac{1}{n}\log|x M_{0}
M_{-1}\cdots
M_{-n+1}|\no\\
&=&\lim_{n\rto}\frac{1}{n}\mbe(\log\| M_{0} M_{-1}\cdots
M_{-n+1}\|).\no
\end{eqnarray*}

In fact we have the following simple but interesting result about
the above the greatest Lyapunov exponents.
\begin{proposition}\label{lja}
  Suppose that condition  \textbf{(C1)}   holds. Then we have
  $$\gamma_L(\overline{M})=\gamma_L( M)=\overline\gamma_L(\overline{M})=\overline\gamma_L(
  M).$$
\end{proposition}
\proof The equality $\gamma_L(M)=\gamma_L(\overline M)$ follows
directly from the similarity of $\overline{M}_n$ and $ M_n.$ Indeed,
for any $n$ we see from the definition of $\overline{M}_n$ and $
M_n$ that $\overline{M}_n=B^{-1} M_nB.$ Since for any matrices $A$
and $B$ we always have $\|AB\|_c\le \|A\|_c\|B\|_c,$ then
\begin{eqnarray*}
\gamma_L(\overline{M})&=&\lim_{n\rto}\frac{1}{n}\log\|\overline{M}_{n-1}\cdots \overline{M}_0\|_c
=\lim_{n\rto}\frac{1}{n}\log\|B^{-1} M_{n-1}\cdots  M_0B\|_c\\
&\le&\lim_{n\rto}\frac{1}{n}\log(\|B^{-1}\|_c\| M_{n-1}\cdots
 M_0\|_c\|B\|_c)=\gamma_L( M).
\end{eqnarray*}
In the same way we have the inverse inequality $\gamma_L(
M)\le\gamma_L(\overline{M})$ to finish  the proof of the first
equality $\gamma_L(\overline{M})=\gamma_L( M).$ The third equality
follows from the same reason as the first equality.

Next we show that
$\gamma_L(\overline{M})=\overline\gamma_L(\overline{M}).$ Note that
for $1\le l\le L-1$
$$\overline{M}_{0}\overline{M}_{-1}\cdots \overline{M}_{-n+1}e_l^T
=\overline{M}_{0}\overline{M}_{-1}\cdots
\overline{M}_{-n+2}((a_{-n+1}(l)e_1^T+e_{l+1}^T).$$ Then one has
that
\begin{eqnarray}
&&\lim_{n\rto}\frac{1}{n}\mbe(\log e_1\overline{M}_{0}\overline
M_{-1}\cdots \overline{M}_{-n+1}e_l^T)\no\\
  &&\hspace{0.5cm}=\lim_{n\rto}\frac{1}{n}\mbe(\log e_1 \overline{M}_{0}\overline{M}_{-1}\cdots
\overline{M}_{-n+2}((a_{-n+1}(l)e_1^T+e_{l+1}^T)\no\\
&&\hspace{0.5cm}=\max\bigg\{\lim_{n\rto}\frac{1}{n}\mbe(\log
a_{-n+1}(l)e_1 \overline{M}_{0}\overline{M}_{-1}\cdots
\overline{M}_{-n+2}e_1^T),\no\\
&&\hspace{3cm}\lim_{n\rto}\frac{1}{n}\mbe(\log e_1
\overline{M}_{0}\overline{M}_{-1}\cdots
\overline{M}_{-n+2}e_{l+1}^T)\bigg\}\no\\
&&\hspace{0.5cm}\ge \lim_{n\rto}\frac{1}{n}\mbe\left(\log e_1
\overline{M}_{0}\overline{M}_{-1}\cdots
\overline{M}_{-n+2}e_{l+1}^T\right)\no\\
&&\hspace{0.5cm}= \lim_{n\rto}\frac{1}{n}\mbe\left(\log e_1
\overline{M}_{0}\overline{M}_{-1}\cdots
\overline{M}_{-n+1}e_{l+1}^T\right).\no
\end{eqnarray}
Therefore it follows that
\begin{eqnarray}
  \overline\gamma_L(\overline{M})\hspace{-0.2cm}&=&\hspace{-0.2cm}\lim_{n\rto}\frac{1}{n}\mbe(\log
  \|\overline{M}_{0}\overline M_{-1}\cdots \overline{M}_{-n+1}\|)\no\\
\hspace{-0.2cm}&=&\hspace{-0.2cm}\lim_{n\rto}\frac{1}{n}\mbe(\log
  e_1\overline{M}_{0}\overline M_{-1}\cdots \overline{M}_{-n+1}\textbf{1}^T)\no\\
\hspace{-0.2cm}&=&\hspace{-0.2cm}\lim_{n\rto}\frac{1}{n}\mbe\Big(\log\sum_{l=1}^L
  e_1\overline{M}_{0}\overline M_{-1}\cdots
  \overline{M}_{-n+1}e_l^T\Big)\no\\
  \hspace{-0.2cm}&=&\hspace{-0.2cm}\max_{1\le l\le L}\lim_{n\rto}\frac{1}{n}\mbe\Big(\log\sum_{l=1}^L
  e_1\overline{M}_{0}\overline M_{-1}\cdots
  \overline{M}_{-n+1}e_l^T\Big)\no\\
\hspace{-0.2cm}&=&\hspace{-0.2cm}\lim_{n\rto}\frac{1}{n}\mbe\Big(\log
  e_1\overline{M}_{0}\overline M_{-1}\cdots
  \overline{M}_{-n+1}e_1^T\Big).\no
\end{eqnarray}
Then one follows from stationarity that
\begin{equation}\label{equ1}
\overline\gamma_L(\overline{M})=\lim_{n\rto}\frac{1}{n}\mbe\Big(\log
  e_1\overline{M}_{n-1}\cdots
  \overline{M}_{1}\overline{M}_0e_1^T\Big).
\end{equation}
 But on the other hand
\begin{eqnarray}\label{equ2}
\gamma_L(\overline
M)\hspace{-0.2cm}&=&\hspace{-0.2cm}\lim_{n\rto}\frac{1}{n}\mbe(\log
\textbf{1}\overline{M}_{n-1}\cdots
\overline{M}_1\overline{M}_0e_1^T)=\lim_{n\rto}\frac{1}{n}\mbe\z(\log\Big(
\sum_{l=1}^Le_1\overline{M}_{n-l}\cdots \overline{M}_1\overline{M}_0e_1^T\Big)\y)\no\\
\hspace{-0.2cm}&=&\hspace{-0.2cm}\max_{1\le l\le
L}\lim_{n\rto}\frac{1}{n}\mbe\z(\log( e_1\overline{M}_{n-l}\cdots
\overline{M}_1\overline{M}_0e_1^T)\y)\no\\
\hspace{-0.2cm}&=&\hspace{-0.2cm}\lim_{n\rto}\frac{1}{n}\mbe\z(\log(
e_1\overline{M}_{n-1}\cdots \overline{M}_1\overline{M}_0e_1^T)\y).
\end{eqnarray}
Then (\ref{equ1}) and (\ref{equ2}) imply that
$\gamma_L(\overline{M})=\overline{\gamma}_L(\overline{M}).$ \qed

Since $\gamma_L(\overline{M})=\gamma_L(
M)=\overline\gamma_L(\overline{M})=\overline\gamma_L(
  M),$ we
write all of them  as $\gamma_L$ in the remainder of our article. We
recall that $\gamma_L$
 characterizes the transience and recurrence of RWRE $X_n$ (see Theorem A of Section \ref{intr}).

Recall that $Z(-k,-m)$ is the $m-k$-th generation of the branching
process beginning at time $-k.$ Let
\begin{equation}\label{yk}
    Y_{-k}=\sum_{m=k+1}^\infty Z(-k,-m)
\end{equation}
being the total number of progeny of the immigrant at times $-k.$
Sinece for $m>k,$ $E_\omega(Z(-k,-m))=M_{-k}M_{-k-1}\cdots
M_{-m+1},$  one has that
\begin{equation}\label{etak}
    \eta_{-k}:=\sum_{m=k+1}^\infty M_{-k}... M_{-m+1},
\end{equation}
is the expectation matrix of $Y_{-k}.$ In next proposition we find
that the projection of $\eta_{0}x_0$ on different directions
$e_l,1\le l\le L,$ that is $e_l\eta_{0}x_0,$ have the same
distribution up to  certain linear transformations.
\begin{proposition}\label{tran}
  Let $x_0=(2,1,...,1)^T\in \mathbb{R}^L$ and $\overline{x}_0=(2,-1,0,...,0)^T\in \mathbb{R}^L.$ Then we
  have\\
i) for all $2\le l\le L,$
  $$\sum_{n=1}^\infty e_l  M_0M_{-1}\cdots
 M_{-n+1}x_0\overset{\mathscr{D}}{=}l+l\sum_{n=1}^\infty
e_1 M_0M_{-1}\cdots M_{-n+1}x_0;$$\\
ii) for all $2<l\le L$
    \begin{equation}
   \sum_{n=1}^\infty e_l\overline{M}_{0}\overline{M}_{-1}\cdots
\overline{M}_{-n+1}\overline{x}_0\overset{\mathscr{D}}{=}
1+\sum_{n=1}^\infty e_1\overline{M}_{0}\overline{M}_{-1}\cdots
\overline{M}_{-n+1}\overline{x}_0,\no
  \end{equation}
  and
  $$\sum_{n=1}^\infty e_2\overline{M}_{0}\overline{M}_{-1}\cdots
\overline{M}_{-n+1}\overline{x}_0\overset{\mathscr{D}}{=}
2+\sum_{n=1}^\infty e_1\overline{M}_{0}\overline{M}_{-1}\cdots
\overline{M}_{-n+1}\overline{x}_0.$$
 \end{proposition}
\proof We mention that all equalities in distribution in our proof
follow from the stationarity of the environment. For fixed  $1<l\le
L,$
\begin{eqnarray}
\sum_{n=1}^\infty e_l  M_0M_{-1}\cdots  M_{-n+1}x_0
\hspace{-0.2cm}&=&\hspace{-0.2cm}\sum_{n=1}^\infty e_l
B\overline{M}_0\overline{M}_{-1}\cdots
\overline{M}_{-n+1}B^{-1}x_0\no\\
\hspace{-0.2cm}&=&\hspace{-0.2cm}\sum_{n=1}^\infty
\sum_{k=1}^le_k\overline{M}_0\overline{M}_{-1}\cdots
\overline{M}_{-n+1}\overline{x}_0. \no
\end{eqnarray}
 For $l=2,$ note that
\begin{eqnarray}\label{p1}
\sum_{n=1}^\infty e_2\overline{M}_0\overline{M}_{-1}\cdots
\overline{M}_{-n+1}\overline{x}_0\hspace{-0.2cm}&=&\hspace{-0.2cm}2+\sum_{n=1}^\infty
e_1\overline{M}_{-1}\overline{M}_{-2}\cdots \overline{M}_{-n}\overline{x}_0\no\\
\hspace{-0.2cm}&\overset{\mathscr{D}}{=}&\hspace{-0.2cm}2+\sum_{n=1}^\infty
e_1\overline{M}_{0}\overline{M}_{-1}\cdots
\overline{M}_{-n+1}\overline{x}_0
\end{eqnarray}Then we have
\begin{eqnarray}
&&\sum_{n=1}^\infty e_2  M_0M_{-1}\cdots
 M_{-n+1}x_0=\sum_{n=1}^\infty
(e_1+e_2)\overline{M}_{0}\overline{M}_{-1}\cdots
\overline{M}_{-n+1}\overline{x}_0\no\\
&&\hspace{1.5cm}\overset{\mathscr{D}}{=}2+2\sum_{n=1}^\infty
e_1\bpim\overline{x}_0 =2+2\sum_{n=1}^\infty e_1B^{-1} \pim B\overline{x}_0\no\\
 &&\hspace{1.5cm}=2+2\sum_{n=1}^\infty e_1\pim x_0.\no
\end{eqnarray}
For $2<l\le L$ note that
\begin{eqnarray}\label{p2}
\sum_{n=1}^\infty
e_l\bpim\overline{x}_0\hspace{-0.2cm}&=&\hspace{-0.2cm}1+
\sum_{n=1}^\infty e_1\overline{M}_{-l+1}\overline{M}_{-l}\cdots\overline{M}_{-l-n+2}\no\\
\hspace{-0.2cm}&\overset{\mathscr{D}}{=}&\hspace{-0.2cm}1+\sum_{n=1}^\infty
e_1\bpim\overline{x}_0.
\end{eqnarray}Then we have
\begin{eqnarray}
&&\sum_{n=1}^\infty e_l \pim
x_0\overset{\mathscr{D}}{=}l+l\sum_{n=1}^\infty
e_1\bpim\overline{x}_0 \no\\
&&\hspace{1cm}=l+l\sum_{n=1}^\infty e_1B^{-1}\pim
B\overline{x}_0=l+l\sum_{n=1}^\infty e_1 \pim x_0.\no
\end{eqnarray}
Then the first assertion follows. The second assertion was also
proved (see the above equations (\ref{p1}) and (\ref{p2})).\qed
{\center \section{An alternative proof of the LLN of
$X_n$}\label{alln}} For $ n\ge 0$ define
$\overline{\omega}(n)=\theta^{X_n}\omega.$ Then
$\{\overline{\omega}(n)\}$ is a Markov chain with transition kernel
$$\overline{P}(\omega,d\omega')=\omega_0(1)\delta_{\theta\omega=\omega'}+\sum_{l=1}^L\omega_0(-l)\delta_{\theta^{-l}\omega=\omega'}.$$ In \cite{br02} an \textbf{(IM)} condition is
said to be satisfied if there is  $\pi(\omega)$ such that
$$\int \tilde{\pi}(\omega)\mbp(d\omega)=1\mbox{ and }
\tilde{\pi}(\omega)=\overline{P}*\tilde{\pi}(\omega),$$ where
$\tilde\pi(\omega)=\pi(\omega)[\mathbb{E}(\pi(\omega))]^{-1}.$ Under
\textbf{(IM)} condition Br\'{e}mont showed an LLN of $\{X_n\}$ in
\cite{br02}. But the  \textbf{(IM)} condition was not given directly
in the words of environment $\omega.$ So one has to check the
existence of the invariant density $\pi(\omega).$ In \cite{br02},
Br\'{e}mont showed the existence of $\pi(\omega)$  by analyzing its
definition and the transition probability of the walk.

What makes difference in our article is that, with the help of the
branching structure, we specify the invariant density $\pi(\omega)$
directly by analyzing a multitype branching process. Therefore we
can avoid introducing the \textbf{(IM)} condition and show directly
that   $\{X_n\}$ satisfies an LLN with a positive speed under the
assumption $``\mathbb{E}(\pi(\omega))<\infty"$. Also the speed has
 a simple explicit form $[\mathbb{E}(\pi(\omega))]^{-1}.$

The result of this section  was stated in Theorem \ref{lln} in the
introduction section. Before giving the proof, we explain how we get
the explicit expression of the invariant density $\pi(\omega)$ from
the MBPREI  constructed in Section \ref{dec}.

For $i\le 0$ define $N_i=\#\{0\le k\le T_1, X_k=i\}.$ Note that
conditioned on the event $\{X_n\rto\},$ $N_i=U^1_{i,1}+|U^1_{i-1}|.$
Then omitting the superscript ``$1$", for $i< 0,$ we have
\begin{eqnarray}
  E_\omega(N_i\big|U_i,U_{i+1},...,U_0)=U_{i,1}+|U_i M_i|.\no
\end{eqnarray}
Hence
\begin{eqnarray}\label{thes}
  E_\omega(N_i)\hspace{-0.2cm}&=&\hspace{-0.2cm}e_1 M_0\cdots M_{i+1}e_1^T+|e_1 M_0\cdots M_{i+1} M_i|\no\\
\hspace{-0.2cm}&=&\hspace{-0.2cm} e_1 M_0\cdots M_{i+1}e_1^T+e_1
M_0\cdots M_{i+1}
 M_i\textbf{1}^T\no\\
\hspace{-0.2cm}&=&\hspace{-0.2cm} (1+a_i(1))e_1 M_0\cdots
M_{i+1}\textbf{1}^T=\frac{1}{\omega_i(1)}e_1B\overline{M}_0\cdots
\overline{M}_{i+1}B^{-1}\textbf{1}^T\no\\
\hspace{-0.2cm}&=&\hspace{-0.2cm}\frac{1}{\omega_i(1)}e_1\overline{M}_0\cdots
\overline{M}_{i+1}e_1^T.
\end{eqnarray}
Note also that
$E_\omega(N_0)=1+E_\omega(|U_{-1}|)=1+\sum_{l=1}^lb_0(l)=\frac{1}{\omega_0(1)}.$
Then
\begin{eqnarray}
E_\omega(
T_1)\hspace{-0.2cm}&=&\hspace{-0.2cm}E_\omega\z(1+\sum_{i=1}^\infty
U_{-i,1}+|U_{-i}|\y)=\sum_{i=0}^\infty
E_\omega(N_{-i})\no\\
\hspace{-0.2cm}&=&\hspace{-0.2cm}\frac{1}{\omega_0(1)}+\sum_{i=1}^\infty\frac{1}{\omega_i(1)}e_1\overline{M}_0\cdots
\overline{M}_{-i+1}e_1^T.\no
\end{eqnarray}

Then we can define
\begin{equation}\label{pi1}\pi(\omega):=\frac{1}{\omega_0(1)}\z(1+\sum_{i=1}^\infty
e_1\overline{M}_i\cdots \overline{M}_1e_1^T\y).\end{equation}
\begin{remark}
  Indeed, one can simply follows from (\ref{tn}) that
 $T_1=1+\sum_{i< 0}U_ix_0$ implying that
 $$E_\omega(
T_1)=1+\sum_{i<0}E_\omega(U_{i}(2,1,...,1)^T)=1+\sum_{i=1}^\infty
e_1M_0\cdots M_{-i+1}(2,1,...,1)^T.$$ Therefore one can also define
\begin{equation}\label{pi2}\pi(\omega)=1+\sum_{i=1}^\infty
e_1M_i\cdots M_{1}(2,1,...,1)^T.\end{equation} One sees that the
right-hand sides of (\ref{pi1}) and (\ref{pi2}) have different
forms. But it follows from the second line of (\ref{thes}) that they
are the same indeed.
\end{remark}
  \noindent \textbf{Proof of Theorem \ref{lln}:} Since the MBPREI $\{U_i\}$ makes sense only  on the event $\{X_n\rto\}$  we first show that
   $\gamma_L<0,$ implying that $P$-a.s., $X_n\rto.$
  Indeed for any $i>0$ by Jensen's inequality  we have
  \begin{eqnarray*}
\mbe(e_1\overline{M}_i\cdots \overline{M}_{1}e_1^T)=\mbe\z(e^{\log
e_1\overline{M}_i\cdots \overline{M}_{1}e_1^T}\y)\ge e^{\mbe\z(\log
e_1\overline{M}_i\cdots \overline{M}_{1}e_1^T\y)}.
  \end{eqnarray*}
But it follows from (\ref{gm})  that
\begin{eqnarray}
\gamma_L\hspace{-0.2cm}&=&\hspace{-0.2cm}\lim_{n\rto}\frac{1}{n}\mbe(\log
\textbf{1}\overline{M}_{n-1}\cdots
\overline{M}_0e_1^T)=\lim_{n\rto}\frac{1}{n}\mbe\z(\log\Big(
\sum_{l=1}^Le_1\overline{M}_{n-l}\cdots \overline{M}_0e_1^T\Big)\y)\no\\
\hspace{-0.2cm}&=&\hspace{-0.2cm}\max_{1\le l\le
L}\lim_{n\rto}\frac{1}{n}\mbe\z(\log( e_1\overline{M}_{n-l}\cdots
\overline{M}_0e_1^T)\y)=\lim_{n\rto}\frac{1}{n}\mbe\z(\log(
e_1\overline{M}_{n-1}\cdots \overline{M}_0e_1^T)\y).\no
\end{eqnarray}
Therefore as $i\rto,$ \begin{eqnarray}e^{\gamma_L i}\sim
e^{\mbe\z(\log e_1\overline{M}_i\cdots \overline{M}_{1}e_1^T\y)}\le
\mbe(e_1\overline{M}_i\cdots \overline{M}_{1}e_1^T)\rightarrow 0,\no
  \end{eqnarray} since $\mbe(\pi(\omega))<\infty.$
  Then we have $\gamma_L<0.$ (i) is proved.
  We just give the idea of the
  proof (ii) and (iii) since it is similar to Zeitouni \cite{ze04} (see the second version of the proof of Theorem 2.1.9).
   It follows from the stationarity that $E(T_1)=\mbe
   (\pi(\omega))<\infty.$ Define $$Q(B)=E\z(\sum_{i=0}^{T_1-1}1_{\{\overline{\omega}(i)\in B\}}\y), \quad
   \overline{Q}(B)=\frac{Q(B)}{Q(\Omega)}=\frac{Q(B)}{E(T_1)}.$$
  Then $Q(\cdot)$ is invariant under kernel $\overline{P},$ that is
$$Q(B)=\iint 1_{\omega'\in B}\overline{P}(\omega,d\omega')Q(d\omega),$$
and
    $\frac{dQ}{d\mbp}=\sum_{i\le
  0}N_i=\pi(\omega).$ Then (\ref{invd}) is proved. Also under $\overline{Q}\otimes P_\omega$ the
  sequence $\{\overline{\omega}(n)\}$ is stationary and ergodic.
  Define the local drift $d(x,\omega)=E_{x,\omega}(X_1-x).$ Then
  \begin{eqnarray}
X_n\hspace{-0.2cm}&=&\hspace{-0.2cm}\sum_{i=1}^n\z(X_i-X_{i-1}-d(X_{i-1},\omega)\y)+\sum_{i=1}^nd(X_{i-1},\omega)\no\\
\hspace{-0.2cm}&:=&\hspace{-0.2cm}
\overline{R}_n+\sum_{i=1}^nd(X_{i-1},\omega)\no
  \end{eqnarray}
where $\{\overline{R}_n\}$ is a $P_{\omega}$-martingale and
$P$-a.s., $\frac{\overline{R}_n}{n}\rightarrow 0.$ We have from the
ergodicity under $\overline{Q}\otimes P_{\omega}$ that $P$-a.s.,
$$\lim_{n\rto}\frac{X_n}{n}=\lim_{n\rto}\frac{1}{n}\sum_{i=1}^n d(0,\overline{\omega}(i-1))=E_{\overline{Q}}(d(0,\overline{\omega}(0))).$$
But
\begin{eqnarray}\label{drif}
&&E_{\overline{Q}}(d(0,\overline{\omega}(0)))=\frac{1}{\mbe(\pi(\omega))}\mbe\z(\pi(\omega)\Big(\omega_0(1)-\sum_{l=1}^Ll\omega_0(-l)\Big)\y)\no\\
&&\hspace{1cm}=\frac{1}{\mbe(\pi(\omega))}\mbe\z(\frac{1}{\omega_0(1)}\Big(1+\sum_{i=1}^\infty
e_1\overline{M}_i\cdots
\overline{M}_{1}e_1^T\Big)\Big(\omega_0(1)-\sum_{l=1}^Ll\omega_0(-l)\Big)\y)\no\\
&&\hspace{1cm}=\frac{1}{\mbe(\pi(\omega))}\mbe\z(\Big(1+\sum_{i=1}^\infty
e_1\overline{M}_i\cdots \overline{M}_{1}e_1^T\Big)\Big(1-\sum_{l=1}^La_0(l)\Big)\y)\no\\
&&\hspace{1cm}=\frac{1}{\mbe(\pi(\omega))}\mbe\z(1+\sum_{i=1}^\infty
e_1\overline{M}_i\cdots
\overline{M}_{1}e_1^T-\sum_{l=1}^La_0(l)(1+\sum_{i=1}^\infty
e_1\overline{M}_i\cdots \overline{M}_{1}e_1^T)\y).
\end{eqnarray}
Note that
\begin{eqnarray}\label{toz}
\sum_{i=1}^\infty e_1\overline{M}_i\cdots
\overline{M}_{1}e_1^T\hspace{-0.2cm}&=&\hspace{-0.2cm}(a_1(1),...,a_1(L))e_1^T+\sum_{i=1}^\infty(a_{i+1}(1),...,a_{i+1}(L))\overline{M}_{i}\cdots
\overline{M}_1e_1^T\no\\
\hspace{-0.2cm}&\overset{\mathscr{D}}{=}&\hspace{-0.2cm}(a_0(1),...,a_0(L))e_1^T+\sum_{i=1}^\infty(a_0(1),...,a_0(L))\overline{M}_{i}\cdots
\overline{M}_1e_1^T\no\\
\hspace{-0.2cm}&=&\hspace{-0.2cm}
a_0(1)+\sum_{i=1}^\infty\sum_{l=1}^La_0(l)e_l \overline{M}_i\cdots
\overline{M}_1e_1^T\no\\
\hspace{-0.2cm}&\overset{\mathscr{D}}{=}&\hspace{-0.2cm}\sum_{l=1}^L\z(a_0(l)+a_0(l)\sum_{i=1}^\infty
e_1\overline{M}_i\cdots \overline{M}_1e_1^T\y),
\end{eqnarray}
where the last step follows similarly as  (\ref{p1}) and (\ref{p2}).
Substituting (\ref{toz}) to (\ref{drif}), using stationarity we have
$E_{\overline{Q}}(d(0,\overline{\omega}(0)))=\frac{1}{\mbe(\pi(\omega))}.$
\qed


{\center \section{The stable limit law of $X_n$}\label{sta}}

To begin with, we introduce some random variables relate to the
process $\{Z_{-n}\}_{n\ge0}$. First recall that
$Y_{-k}=\sum_{m=k+1}^\infty Z(-k,-m) $ is the total number of
progeny of the immigrant at times $-k,$ and
$\eta_k:=\sum_{m=k+1}^\infty M_{-k}... M_{-m+1}$ is the
corresponding
 expectation random matrix.

Next, let $\nu_0\equiv0,$ and define recursively
\begin{equation}
    \nu_n=\min\{m>\nu_{n-1}:Z_{-m}=\textbf{0}\}\mbox{ for } n>0,\no
\end{equation} being the successive regeneration times of MBPREI $\{Z_{-n}\}_{n\ge0}.$
 For simplicity we write
$\nu_1$ as $\nu.$

Define also
\begin{equation}
    W=\sum_{k=0}^{\nu-1}Z_{-k},\no
\end{equation}
the total number of offspring born before regenerating time $\nu.$

Finally, for $A>0,$ we introduce the stopping time
\begin{equation}
    \sigma=\sigma(A)=\inf\{m:|Z_{-m}|>A\}\no
\end{equation}
which is the time the number of particles of the process $\{Z_n\}$
exceeding $A.$


\subsection{The tail of the expectation of the total number of $\{Z(0,-k)\}_{k\ge 0}$}
To study the limit law of RWRE with bounded jumps $\{X_n\},$ a key
step is to prove that random variable $Wx$ belongs to the domain of
attraction of some $\ka$-stable law. For this purpose it is crucial
to show first that $x\eta_0x_0$
 belongs to the domain of attraction of a $\ka$-stable law for any positive $x\in
 \mathbb{R}^L$. Indeed, we have
\begin{theorem}\label{ep}
  Suppose that $\gamma_L<0.$ Then under Condition C,
  for  $\kappa$ of (\ref{ka}) (see also (\ref{wk}) below) and for some $K_2=K_2(x_0)\in (0,\infty),$
  we have
  \begin{equation}
 \lim_{t\rto} t^\kappa\mbp(x\eta_0x_0\ge t)=K_2|xB|^{\kappa}\no
  \end{equation}
  for all $x\in \mathbb{R}^L$ with positive components such that $|x|>0.$
 \end{theorem}

To proof Theorem \ref{ep},
 we need some classical results of random
matrices in Kesten's paper \cite{kes73}. We rewrite them in terms of
$\{ M_{-n}\}_{n\ge 0}.$ Recall that
$$\gamma_L=\lim_{n\rto}\frac{1}{n}\mbe(\log\| \pim\|)$$
is the greatest Lyapunov exponent of $\{M_i\}_{i\le 0}.$
\begin{theorem}[Kesten \cite{kes73}]\label{kes} Suppose that Condition C holds and $\gamma_L<0.$  Then
\begin{itemize}
  \item[1)] for every $\alpha\in[0,\kappa_0]$ the limits
  \begin{eqnarray}\label{wk1}
    \log\rho(\alpha)\hspace{-0.3cm}&:=&\hspace{-0.3cm}\lim_{n\rto}\frac{1}{n}\log\mbe\z(\parallel \pim\parallel^\alpha\y)
     \end{eqnarray} exist and $\log\rho(\alpha)$ is a strictly convex
  function of $\alpha$. Hence
  \item[2)] there exists a unique $\kappa\in(0,\kappa_0],$ such that
\begin{eqnarray}\label{wk}
  \log\rho(\kappa)\hspace{-0.3cm}&=&\hspace{-0.3cm}\lim_{n\rto}\frac{1}{n}\log\mbe\z(\parallel \pim\parallel^{\kappa}\y)=0.
\end{eqnarray}
\item[3)] Let $\{Q_{-n}\}_{n\ge0},$ with the law of $Q_0$ being $\mathbb{Q},$
be a random sequence of  L-(column) vectors such that $\{
M_{-n},Q_{-n}\}_{n\ge0}$ are i.i.d.. Assume also
$\mathbb{Q}(Q_0=\textbf{0})<1,$ $\mathbb{Q}(Q_0\ge \textbf{0})=1,$
$E_\mathbb{Q}{|Q_0|^\kappa}<\infty$  for $\kappa$ of (\ref{wk}),
where $Q_0\ge \textbf{0}$ means that all components of $Q_0$ are
nonnegative. Then for each $x\in S_{L-1},$ with an abuse use of
notation $\mathbb{P},$ the limit
$$\lim_{t\rto}t^\kappa \mathbb{P}\Big(\sum_{n=1}^\infty
 x\pim Q_{-n}>t\Big)$$ exists and is finite. In particular
there exist constants $K_1=K_1(\mathbb{P}, M,\mathbb{Q})\in
(0,\infty)$ and  $r=r(x, M)\in(0,\infty)$ such that
\begin{eqnarray}\label{li}
  \lim_{t\rto}t^\kappa \mathbb{P}\Big(\sum_{n=1}^\infty
 x\pim Q_{-n}\ge t\Big)= K_1(\mathbb{P}, M,\mathbb{Q})r(x,
M)
\end{eqnarray}
for $x\in S_+.$
\end{itemize}
\end{theorem}
\begin{remark}\label{rk}
  (i) The first two parts of the theorem can be concluded from the
  proof of Theorem 3 (see step 4) of Kesten \cite{kes73}.  The third part
  corresponds to Theorem 4 of Kesten \cite{kes73}.\\
  (ii) We mention that  $|x|$ denotes
$\Big(\sum_{i=1}^Lx_i^2\Big)^\frac{1}{2}$ in \cite{kes73}. But all
proofs go through under $l_1$-norm $|x|:=\sum_{i=1}^L|x_i|.$
\end{remark}
Now we are ready to present the proof of Theorem \ref{ep}.

\noindent{\bf Proof of Theorem \ref{ep}:} Fix $x=(x_1,...,x_L)\in
\mathbb{R}^L$ such that $x_i\ge 0,\ i=1,...,L$
 and $|x|>0.$ We have
 \begin{eqnarray}
   \mbp(x\eta_0x_0\ge t)\hspace{-0.3cm}&=&\hspace{-0.3cm}\mbp\Big(\sum_{n=1}^\infty x M_0M_{-1}\cdots M_{-n+1}x_0\ge
   t\Big)=\mbp\Big(\sum_{n=1}^\infty \sum_{l=1}^Lx_le_l M_0M_{-1}\cdots M_{-n+1}x_0\ge
   t\Big).\no
\end{eqnarray}
 It follows from Proposition \ref{tran} that the rightmost-hand side of
 above expression equals to
 \begin{eqnarray}
  && \mbp\Big(\sum_{l=2}^{L}lx_l+\sum_{l=1}^Llx_l\sum_{n=1}^\infty e_1 \pim x_0\ge
  t\Big)\no\\
  &&\hspace{1.5cm}=\mbp\Big(\sum_{l=2}^{L}lx_l+|xB|\sum_{n=1}^\infty e_1 \pim x_0\ge
  t\Big).\no
 \end{eqnarray}
 Then one gets from the third part of Theorem \ref{kes} that
 \begin{eqnarray}
 &&\lim_{t\rto} t^\kappa\mbp(x\eta_0x_0>t)=\lim_{t\rto} t^\kappa\mbp\Big(|xB|\sum_{n=1}^\infty e_1 \pim x_0\ge
  t\Big)\no\\
 &&\hspace{2cm}=|xB|^\kappa
 K_1(\mbp, M,\delta_{x_0})r(e_1, M)=:K_2|xB|^\kappa,\no
  \end{eqnarray}
which finishes the proof of Theorem \ref{ep}.\qed
  \subsection{The tail of the population size of MBPREI before regeneration}\label{pz}
Recall that $\nu$ is the regeneration time of the MBPREI and
$W=\sum_{0<n\le\nu-1}Z_{-n}$ is the total number of particles born
to the immigrants entering before time $\nu.$ The main purpose of
this section is to find how large the population size will be before
regeneration.
\begin{theorem}\label{siz}
Suppose that   Condition C holds and $\gamma_L<0.$  If  $\kappa> 2,$
then $ E((Wx_0)^2)<\infty;$ if $\kappa\le 2,$ then there exists some
$0<K_3<\infty$ such that
\begin{equation}\label{ws}
\lim_{t\rto}t^\kappa P(Wx_0\ge t)=K_3.
\end{equation}
\end{theorem}
To prove the theorem, we need some preparations. To begin with we
 show that the tail probability of $\nu$ vanishes with
exponential rate. This can follow from Theorem 4.2 of Key
\cite{key87}. But for our MBPREI  condition (iii):
$\mbp(P_\omega(Z(0,-1)=\textbf{0}|Z(0,0)=e_l)>0\ for\
l=1,2,...,L)>0,$ in that theorem does not hold. However, we can
prove the results directly, since we have
\begin{lemma}\label{po}
Suppose that condition \textbf{(C1)} holds and $\gamma_L<0.$
   Then
  $\lim_{m\rto}P(Z_{-m}=v)=\pi(v),$  a probability distribution on
  $\mathbb{Z}^L$ with
 $\pi(\textbf{0})>0.$
\end{lemma}
\proof The first assertion follows as Theorem 3.3 in \cite{key87}.
To show that $\pi(\textbf{0})>0$ we proceed by contradiction. If
$\pi(\textbf{0})=0,$ then $\lim_{m\rto}P(Z_{-m}=\textbf{0})=0.$
Hence we have $\mbp$-a.s.,
\begin{equation}\label{zm1}
\lim_{m\rto}P_\omega(Z_{-m}=\textbf{0})=0.
\end{equation}
For $v\in \mathbb{Z}^L,$ $m>0,$ let
\[q(m,v):=P_\omega(Z_{-m-L}=\textbf{0}|Z_{-m}=v).\]
Then for each $v\in  \mathbb{Z}^L $
\begin{equation}\label{zm2}
    P_\omega(Z_{-m-L}=\textbf{0})\ge
    P_{\omega}(Z_{-m}=v)q(m,v).
\end{equation}
Taken together, (\ref{zm1}) and  (\ref{zm2}) imply that for all
$v\in\mathbb{Z}^L,$ $\mbp$-a.s.,
\begin{equation*}
\lim_{m\rto}P_{\omega}(Z_{-m}=v)q(m,v)=0.
\end{equation*}
It follows by stationarity that $P_{\omega}(Z_{-m}=v)q(m,v)$ and
$P_{\omega}(Z'_{-m}=v)q(0,v)$ have the same distribution, where
$Z'_{-m}=\sum_{k=0}^{m-1}Z(k,0).$ Therefore for all
$v\in\mathbb{Z}^L,$
 $\mbp$-a.s.,
\begin{equation*}
\lim_{m\rto}P_{\omega}(Z'_{-m}=v)q(0,v)=0.
\end{equation*}
Then on the event $\{q(0,v)>0\},$
\begin{equation*}
\lim_{m\rto}P_{\omega}(Z'_{-m}=v)=0.
\end{equation*}
If we can show that
\begin{equation}\label{z0l}
\mbp(q(0,v)>0)=1
\end{equation} then it follows that $\mbp$-a.s.,
$$\sum_{v\in\mathbb{Z}^L}\lim _{m\rto}P_\omega(Z'_{-m}=v)=0$$ which
will contradict that $\pi$ is a probability distribution (Here we
mention  that  $Z'_{-m}$ and $Z_{-m}$ have the same limit
distributions, see Lemma 2.1 and Lemma 3.2 of Key \cite{key87}). It
remains to show (\ref{z0l}). Since
$$\max_{1\le l\le
L}\mbe\Big(\log^+\frac{\omega_0(-l)}{\omega_{0}(1)}\Big)<\infty,$$
$\mbp(\omega_0(1)>0)=1.$ Then we have $\mbp$-a.s.,
$$P_\omega(Z_{-L}=\textbf{0}|Z_0=(v_1,v_2,...,v_L))\ge\prod_{k=0}^{L-1}\omega_{-k}(1)^{1+\sum_{l=k+1}^Lv_l}>0$$
which proves (\ref{z0l}).\qed
\begin{remark}
 We adopted the idea of the proof of Theorem 3.3 of Key \cite{key87} to prove this lemma.
 The only difference here is that  we replace condition $$(iii):
  \mbp(P_\omega(Z(0,-1)=\textbf{0}\big|Z(0,0)=e_l)>0\ for\
  l=1,2,...,L)>0,$$
  in that theorem with $\mbp(\omega_0(1)>0)>0,$ which is implied in
  condition \textbf{(C1)}.
\end{remark}
With Lemma \ref{po} in hands, the next theorem follows verbatim as
Theorem 4.2 in \cite{key87}.
\begin{theorem}[Key\cite{key87}]\label{nu}
  Suppose that $\gamma_L<0$ and that condition \textbf{(C3)} holds. Then there
  exist positive constants $K_4$ and $K_5$ such that
  \begin{equation}
    P(\nu>t)<K_4\exp(-K_5t).\no
  \end{equation}
\end{theorem}
In the following three lemmas, i.e., Lemma \ref{ys}, Lemma \ref{zsk}
and Lemma \ref{sz}, some estimations for the related probabilities
and moments are given. The proofs are technical and follow almost
verbatim as
 Kesten-Kozlov-Spitzer \cite{kks75}. But the proofs will be  long journeys.
 Therefore, for continuity consideration, we delay the proofs of these lemmas to
 the Appendix section.

 In  Lemma \ref{ys}, Lemma \ref{zsk}
and Lemma \ref{sz}, we always make the assumption that all
conditions of Theorem \ref{siz} hold.
\begin{lemma}\label{ys}
  If $\kappa\le2,$ then there exists for all $\epsilon>0$ an
  $A_0=A_0(\epsilon)<\infty$ such that
  \begin{equation}
    P\bigg(\sum_{\sigma\le k<\nu}|Y_{-k}|\ge \epsilon x\bigg)\le \epsilon
    x^{-\kappa} \mbox{    for } A\ge A_0(\epsilon).\no
  \end{equation}
  \end{lemma}

\begin{lemma}\label{zsk}
  If $\kappa\le2,$ then for fixed $A$
  \begin{equation}\label{zsi}
    E\z(|Z_{-\sigma}|^{\kappa};\sigma<\nu\y)<\infty.
  \end{equation}
  If $\kappa>2$ then
  \begin{equation}
    E(|W|^2)<\infty.\no
  \end{equation}
  \end{lemma}
Next we introduce $$S_{-\sigma,-m}=\mbox{number of progeny alive at
time }-m\mbox{ of the $Z_{-\sigma}$ particles present at
$-\sigma>-m$.}$$ Let $S_{-\sigma,-\sigma}=Z_{-\sigma},$  and
$$S_{-\sigma}=\sum_{m=\sigma}^\infty S_{-\sigma,-m}=Z_{-\sigma}+\mbox{total progeny of the $Z_{-\sigma}$ particles at $-\sigma.$}$$
\begin{lemma}\label{sz}
  If $\ka\le 2,$ then there exists for all $\epsilon>0$ an
  $A_1=A_1(\epsilon)$ such that  for $A>A_1$
  \[P\Big(\Big|\sum_{m=\sigma}^\infty\Big(S_{-\sigma,-m}-Z_{-\sigma}\prod_{i=\sigma}^{m-1} M_{-i}\Big)\Big|\ge \epsilon x,\sigma<\nu\Big)\le
  \epsilon x^{-\ka}E(|Z_{-\sigma}|^\ka;\sigma<\nu).\]
\end{lemma}

\noindent \textit{Proof of Theorem \ref{siz}}: Since in Lemma
\ref{zsk} we have shown that $ E(|W|^2)<\infty$ when $\ka>2,$ it
follows immediately that $ E((Wx_0)^2)< 4E(|W|^2)<\infty.$ The first
part of the theorem follows. To prove the second part, recall that
$W$ is the number of particles born before
  $-\nu.$ Then on the event $\{\sigma<\nu\}$ we have \[W=\sum_{s=0}^{\sigma-1}Z_{-s}+S_{-\sigma}+\sum_{\sigma\le s<\nu}Y_{-s}.\]
As an immediate corollary of Theorem \ref{nu} we have for all
$\epsilon>0,\ A>0$
\begin{equation}\label{wan1}
  P(W x_0\ge \epsilon x,\sigma(A)\ge\nu)\le P(2A\nu\ge\epsilon
  x)=o(x^{-\ka}),\quad\quad x\rto,
\end{equation}
since $W x_0<|2W|=|2\sum_{t=0}^{\nu-1}Z_{-t}|\le 2A\nu$ on the event
$\{\sigma\ge \nu\}.$ Similarly we have
\begin{equation}\label{wan2}
  P\Big(\Big|\sum_{t=0}^{\sigma-1}Z_{-t} x_0\Big|\ge \epsilon x,\sigma(A)<\nu\Big)\le P(2A\nu\ge\epsilon
  x)=o(x^{-\ka}),\quad\quad x\rto.
\end{equation}
Taken together (\ref{wan1}), (\ref{wan2}) and Lemma \ref{ys} imply
that for sufficiently large $A$ and $x$
\begin{eqnarray}\label{sws}
 && P(\sigma<\nu,S_{-\sigma} x_0\ge x)\no\\
 &&\quad\quad\le P(W x_0\ge x)=P(W x_0\ge x,\sigma<\nu)+P(W x_0\ge x,\sigma\ge\nu)\no\\
 &&\quad\quad\le P(\sigma<\nu,S_{-\sigma} x_0\ge x-\sum_{s=0}^{\sigma-1}Z_{-s} x_0-\sum_{\sigma\le
 s<\nu}Y_{-s} x_0)+\epsilon x^{-\ka}\no\\
 &&\quad\quad\le P(\sigma<\nu,S_{-\sigma} x_0\ge x(1-2\epsilon))+3\epsilon x^{-\ka}.
\end{eqnarray}
Since \begin{eqnarray}
&&P\Big(\sigma<\nu,Z_{-\sigma}\sum_{t=\sigma}^\infty\prod_{i=\sigma}^{t-1}
M_{-i} x_0\ge
(1+\epsilon)x\Big)\no\\
&&\quad\quad\le P\Big(\sigma<\nu, \sum_{t=\sigma}^\infty
S_{-\sigma,-t} x_0\ge
x\Big)+P\Big(\sigma<\nu,\Big|\Big(Z_{-\sigma}\sum_{t=\sigma}^\infty\prod_{i=\sigma}^{t-1}
M_{-i}-\sum_{t=\sigma}^\infty
S_{-\sigma,-t}\Big) x_0\Big|\ge\epsilon x\Big)\no\\
&&\quad\quad\le P\Big(\sigma<\nu, \sum_{t=\sigma}^\infty
S_{-\sigma,-t} x_0\ge
x\Big)+P\Big(\sigma<\nu,2\Big|Z_{-\sigma}\sum_{t=\sigma}^\infty\prod_{i=\sigma}^{t-1}
M_{-i}-\sum_{t=\sigma}^\infty S_{-\sigma,-t}\Big|\ge\epsilon
x\Big),\no
\end{eqnarray}
and \begin{eqnarray} &&P(\sigma<\nu,S_{-\sigma}
x_0>x(1-2\epsilon))\le P\Big(\sigma<\nu,
Z_{-\sigma}\sum_{t=\sigma}^\infty\prod_{i=\sigma}^{t-1} M_{-i}
x_0\ge
(1-3\epsilon)x \Big)\no\\
&&\quad\quad\quad\quad\quad\quad
+P\Big(\sigma<\nu,\Big|\Big(Z_{-\sigma}\sum_{t=\sigma}^\infty\prod_{i=\sigma}^{t-1}
M_{-i}-\sum_{t=\sigma}^\infty
S_{-\sigma,-t}\Big) x_0\Big|\ge\epsilon x\Big)\no\\
&&\quad\quad\le P\Big(\sigma<\nu,
Z_{-\sigma}\sum_{t=\sigma}^\infty\prod_{i=\sigma}^{t-1} M_{-i}
x_0\ge
(1-3\epsilon)x \Big)\no\\
&&\quad\quad\quad\quad\quad\quad
+P\Big(\sigma<\nu,2\Big|Z_{-\sigma}\sum_{t=\sigma}^\infty\prod_{i=\sigma}^{t-1}
M_{-i}-\sum_{t=\sigma}^\infty S_{-\sigma,-t}\Big|\ge\epsilon
x\Big),\no
\end{eqnarray}
  then it follows from Lemma \ref{sz} and (\ref{sws}) that for
  sufficiently large $A$
\begin{eqnarray}\label{sws1}
&&P\Big(\sigma<\nu,Z_{-\sigma}\sum_{t=\sigma}^\infty\prod_{i=\sigma}^{t-1}
M_{-i} x_0\ge
(1+\epsilon)x\Big)-\epsilon x^{-\ka}E(|Z_{-\sigma}|^\ka;\sigma<\nu)\no\\
&&\quad\quad\le P(W x_0\ge x)\no\\
&&\quad\quad\le P\Big(\sigma<\nu,
Z_{-\sigma}\sum_{t=\sigma}^\infty\prod_{i=\sigma}^{t-1} M_{-i}
x_0\ge (1-3\epsilon)x \Big)+ \epsilon
x^{-\ka}(3+E(|Z_{-\sigma}|^\ka;\sigma<\nu)).
\end{eqnarray}
Since
$$\sum_{t=\sigma}^\infty\prod_{i=\sigma}^{t-1} M_{-i}=I+\eta_{-\sigma},$$
we can write (\ref{sws1}) as
\begin{eqnarray}
&&P(\sigma<\nu,Z_{-\sigma} (I+\eta_{-\sigma}) x_0\ge
(1+\epsilon)x)-\epsilon x^{-\ka}E(|Z_{-\sigma}|^\ka;\sigma<\nu)\no\\
&&\quad\quad\le P(Wx_0\ge x)\no\\
&&\quad\quad\le P(\sigma<\nu, Z_{-\sigma}(I+\eta_{-\sigma}) x_0\ge
(1-3\epsilon)x)+ \epsilon
x^{-\ka}(3+E(|Z_{-\sigma}|^\ka;\sigma<\nu)).\no
\end{eqnarray}
Then, to prove (\ref{ws}), it suffices to prove that for each fixed
$A$ \begin{equation}\label{usoft} 0<\lim_{x\rto}x^\ka
P(\sigma=\sigma(A)<\nu,Z_{-\sigma}(I+\eta_{-\sigma}) x_0\ge
x)=K_2E(|Z_{-\sigma} B|^\ka;\sigma<\nu)<\infty.
\end{equation}
 Indeed, this follows immediately from Lemma
\ref{zsk} and Theorem \ref{ep} since
\begin{eqnarray}
  && \lim_{x\rto}x^\ka
P(\sigma<\nu,Z_{-\sigma}(I+\eta_{-\sigma}) x_0\ge x)\no\\
&&\quad\quad=\lim_{x\rto}x^\ka \int_{|s|\ge A}^\infty
P(\sigma<\nu,Z_{-\sigma}\in ds)P(s(I+\eta_{-\sigma}) x_0\ge
x)\no\\
&&\quad\quad=K_2 \int_{|s|\ge A}^\infty P(\sigma<\nu,Z_{-\sigma}\in
ds)|sB|^\ka\no\\
&&\quad\quad=K_2E(|Z_{-\sigma} B|^\ka;\sigma<\nu)<L^\ka
K_2E(|Z_{-\sigma} |^\ka;\sigma<\nu)<\infty.\no
\end{eqnarray}

Also, we have $$ E(|Z_{-\sigma} B|^\ka;\sigma<\nu)\ge E(|Z_{-\sigma}
|^\ka;\sigma<\nu)\ge
A^{\ka}\mbe(P_{\omega}(|Z_{-1}|>A))>A^{\ka}\mbe(\omega_0(-1)^{A+1}\omega_0(1))>0.$$
Thus Theorem \ref{siz} follows.\qed.

\noindent{\bf Proof of Theorem \ref{main}:} From here on the proof
of Theorem \ref{main} is standard. Recall that  $$
T_n=n+\sum_{i=-\infty}^{n-1}|U^n_i|+\sum_{i=-\infty}^{n-1}U^n_{i,1}=n+\sum_{i=-\infty}^{n-1}U^n_{i}x_0.$$
 When the walk
is transient to the right, it only takes finite steps in
$(-\infty,0),$ i.e., $P$-a.s., $$\sum_{i<0}U_i^nx_0<\infty.$$
Therefore to determine the limit distribution of $T_n,$ we need only
to consider $$ n+\sum_{i=0}^{n-1}U^n_ix_0,$$ which, by Theorem
\ref{edis}, has the same distribution with
$$n+\sum_{t=0}^{n-1}Z_{-t}x_0.$$

In Theorem \ref{siz} we have proved that if $\ka>2,$
$E((Wx_0)^2)<\infty$ while for $\ka\le2,$ $P(W x_0>x)\sim K_3x^\ka,$
as $x\rto.$ Then it follows that (see, Feller \cite{fe71}, Chap.
XVII, Sec. 5 Theorem 2 or Gendenko-Kolmogorov \cite{gk54} Chap.7,
Sec. 35, Theorem 2) $W x_0$ belongs to the domain of attraction of a
$\ka$-stable law. Recall that $\nu_0=0,\nu_1,\nu_2,...$ are the
successive regeneration times of the MBPREI. Now put
$$W_{k}=\sum_{\nu_k\le t<\nu_{k+1}}Z_{-t},$$ then the pairs
$\{(\nu_{k+1}-\nu_k),W_{k}\}_{k\ge0}$ are independent, all with
distribution of the $(\nu, W)$ because $(\nu, W)$ coincides with
$(\nu_1-\nu_0,W_0).$ Then the proof of the theorem is standard. It
follows exactly with the proof of the theorem in
Kesten-Kozlov-Spitzer \cite{kes73}. We will not repeat it here.\qed

{\section*{Appendix: Proofs of Lemma \ref{ys}, Lemma \ref{zsk} and
Lemma \ref{sz}}}

\noindent{\bf Proof of Lemma \ref{ys}:}  Since that
$\sum_{k=1}^\infty k^{-2}=\frac{\pi^2}{6},$ we have
\begin{eqnarray}\label{duli}
P\bigg(\sum_{\sigma\le k<\nu}|Y_{-k}|\ge \epsilon
x\bigg)\hspace{-0.2cm}&=&\hspace{-0.2cm}P\z(\sum_{k=1}^\infty1_{[\sigma\le
k<\nu]}|Y_{-k}|\ge 6\pi^{-2}\epsilon x\sum_{k=1}^\infty k^{-2}\y)\no\\
\hspace{-0.2cm}&\le&\hspace{-0.2cm} \sum_{k=1}^\infty P\z(\sigma\le
k<\nu,|Y_{-k}|\ge \frac{1}{2}\epsilon x k^{-2}\y).
\end{eqnarray}
Note that $Y_{-k}$ is
$\sigma\{\omega_{-k},\omega_{-k-1},\cdots\}$-measurable, the event
$\{\sigma\le k<\nu\}$  being defined in terms of
$Z_0,Z_{-1}\cdots,Z_{-k}$ is
$\sigma\{\omega_0,\omega_{-1}\cdots,\omega_{-k+1}\}$-measurable and
$Y_{-k}$ has the same distribution as $Y_0.$ Then
\begin{eqnarray*}
P\Big(\sum_{\sigma\le k<\nu}|Y_{-k}|\ge \epsilon
x\Big)\hspace{-0.2cm}&\le&\hspace{-0.2cm} \sum_{k=1}^\infty
P\z(\sigma\le k<\nu\y)P\Big(|Y_{-k}|\ge \frac{1}{2}\epsilon x
k^{-2}\Big)\\
\hspace{-0.2cm}&=&\hspace{-0.2cm} \sum_{k=1}^\infty P\z(\sigma\le
k<\nu\y)P\Big(|Y_0|\ge \frac{1}{2}\epsilon x k^{-2}\Big).
\end{eqnarray*}
Thus if we can prove that
\begin{equation}\label{y0}
    P(|Y_0|\ge x)\le K_6x^{-\kappa}
\end{equation} for some $K_6<\infty,$ then it follows that
\begin{equation}
\begin{split}
P\bigg(\sum_{\sigma\le k<\nu}|Y_{-k}|\ge \epsilon x\bigg)&\le
x^{-\kappa}2^{\kappa}\epsilon^{-\kappa}K_6\sum_{k=1}^\infty
k^{2\kappa}P\z(\sigma\le k<\nu\y)\\
&\le
x^{-\kappa}2^\ka\epsilon^{-\kappa}K_6E\z(\nu^{2\kappa+1};\sigma<\nu\y)\le
\epsilon x^{-\kappa}\no
\end{split}
\end{equation}
for $A$ large enough since $E\z(\nu^{2\kappa+1}\y)<\infty $ by Lemma
\ref{ys} and $\sigma(A)\rightarrow \infty$ in probability as
$A\rightarrow\infty.$

To prove (\ref{y0}), observe that $\eta_{-m}= M_{-m}(I+\eta_{-m-1})$
and consequently with $Z(0,0)=\textbf{0}$
\begin{equation}
Y_0=\sum_{m=1}^\infty Z(0,-m)=\sum_{m=1}^\infty\z(Z(0,-m)-Z(0,-m+1)
M_{-m+1}\y)(I+\eta_{-m}).\no
\end{equation}
Let $e_0=(\frac{1}{L},...,\frac{1}{L}).$ Using the independence of
$(I+\eta_{-m})$ and $ M_{-m+1},Z(0,-m+1),Z(0,-m),$ similarly as
(\ref{duli}) we have
\begin{eqnarray}
P(|Y_0|\ge x)&\le&\sum_{m=1}^\infty
P\z(|\z(Z(0,-m)-Z(0,-m+1) M_{-m+1}\y)(I+\eta_{-m})|\ge\frac{1}{2}m^{-2}x\y)\no\\
&\le&\sum_{m=1}^\infty\int P\z(|Z(0,-m)-Z(0,-m+1) M_{-m+1}|\in
ds\y)\no\\
&\times&P\z(e_0( I+\eta_0)e_0^T\ge(2sL^2m^2)^{-1}x\y).\no
\end{eqnarray}
From Theorem \ref{kes} there exists a $0<K_7<\infty$ for which
\begin{equation}\label{k7}
P\z(e_0( I+\eta_0)e_0^T\ge(2sL^2m^2)^{-1}x\y)\le
K_7(2sL^2m^2)^{\kappa}x^{-\kappa}.
\end{equation} Then it follows that
\begin{eqnarray}\label{seri}
&&P(|Y_0|\ge x)\le
x^{-\kappa}2^{\kappa}L^{2\kappa}K_7\sum_{m=1}^\infty
m^{2\kappa}E\z(|Z(0,-m)-Z(0,-m+1) M_{-m+1}|^{\kappa}\y)\no\\
&&\quad\le x^{-\kappa}2^{\kappa}L^{2\kappa}K_7\sum_{m=1}^\infty
m^{2\kappa}\mbe\z(E_\omega\z(|Z(0,-m)-Z(0,-m+1)
M_{-m+1}|^2\y)^{\frac{\kappa}{2}}\y)
\end{eqnarray}
 since by assumption of the lemma $\kappa\le2.$ We prove next the
convergence of the last series in above expression to complete the
proof of the lemma. For this purpose note that
\begin{eqnarray}
&&E_\omega\z(|Z(0,-m)-Z(0,-m+1) M_{-m+1}|^2|Z(0,-m+1)\y)\no\\
&&\quad\quad=|Z(0,-m+1)|\Big(\sum_{l=1}^L(b_{-m+1}(l)+b_{-m+1}^2(l))+2\sum_{1\le
l<k\le L}b_{-m+1}(l)b_{-m+1}(k)\Big)\no\\
&&\quad\quad=:|Z(0,-m+1)|R( M_{-m+1}). \label{fc1}
\end{eqnarray}
Then it follows that
\begin{eqnarray}
&&\mbe\z(E_\omega\z(|Z(0,-m)-Z(0,-m+1)
M_{-m+1}|^2\y)^{\frac{\kappa}{2}}\y)\no\\
&&\hspace{1cm}=\mbe(\{E_\omega(|Z(0,-m+1)|
R( M_{-m+1})\}^{\frac{\kappa}{2}}))\no\\
&&\hspace{1cm}=\mbe\Big(|e_1 M_0M_{-1}...
M_{-m+2}|^{\frac{\kappa}{2}} R(
M_{-m+1})^{\frac{\kappa}{2}}\Big)\no\\
&&\hspace{1cm}=\mbe\z(|e_1 M_0M_{-1}...
M_{-m+2}|^{\frac{\kappa}{2}}\y)\mbe(R( M_0)^{\frac{\kappa}{2}}),\no
\end{eqnarray}using independence and stationarity in the last step.
Since $\ka<2,$ then condition \textbf{(C3)} implies that $\mbe(R(
M_0)^\frac{\ka}{2})<\infty,$ and  1), 2) of Theorem \ref{kes} imply
that
\begin{eqnarray*}
&&\mbe\Big(|e_1 M_0M_{-1}\cdots
 M_{-m+2}|^{\frac{\ka}{2}}\Big)
\le L^\frac{\ka}{2}\mbe\Big(\| M_0M_{-1}\cdots
 M_{-m+2}\|^{\frac{\ka}{2}}\Big)\sim L^\frac{\ka}{2}e^{-c(m-1)}
\end{eqnarray*}
 as $m$ tends to $\infty$ for some
constant $c>0.$ Thus the convergence of the last series in
(\ref{seri}) follows.\qed

\noindent{\bf Proof of Lemma \ref{zsk}:} We have on $\{\sigma<\nu\}$
  \begin{equation}\label{zzz}
|Z_{-\sigma}|=(|Z_{-\sigma+1}|+1)\frac{|Z_{-\sigma}|}{|Z_{-\sigma+1}|+1}\le
(A+1)\frac{|Z_{-\sigma}|}{|Z_{-\sigma+1}|+1}
  \le(A+1)\sum_{1\le m\le \nu}\frac{|Z_{-m}|}{|Z_{-m+1}|+1}.
  \end{equation}
  As a matter of fact, if $\kappa\ge 1$
  \begin{eqnarray}\label{zss}
&&\z(E\z(|Z_{-\sigma}|^\ka;\sigma<\nu\y)\y)^{\frac{1}{\ka}}\le
(A+1)\Big(E\bigg[\bigg(\sum_{m\ge1}\frac{|Z_{-m}|}{|Z_{-m+1}|+1}1_{[m\le\nu]}\bigg)^\ka\bigg]\Big)^{\frac{1}{\ka}}\no\\
&&\quad\quad\le
(A+1)\sum_{m\ge1}\Big(E\bigg(\Big(\frac{|Z_{-m}|}{|Z_{-m+1}|+1}\Big)^\ka1_{[m\le\nu]}\bigg)\Big)^{\frac{1}{\ka}}.
  \end{eqnarray}
Conditioned on $\omega$ and $Z_{-m+1}$ we have
$$|Z_{-m}|\le \sum_{j=1}^{|Z_{-m+1}|+1}(|V_j|+1),$$
where

$$P_\omega(V_j=(u_1,\cdots,u_L))
=\frac{(u_1+\cdots+u_L)!}{u_1!\cdots
u_L!}\omega_{-m+1}(-1)^{u_1}\cdots\omega_{-m+1}(-L)^{u_L}\omega_{-m+1}(1)$$
 and $V_j,\ j=1,2,...$ are i.i.d..
Then for $1\le \ka\le 2$ we have
\begin{eqnarray}\label{zmm}
  &&(E_\omega(|Z_{-m}|^\ka\big|Z_0,Z_{-1},\cdots,Z_{-m+1}))^{\frac{1}{\ka}}\le
  \sum_{j=1}^{|Z_{-m+1}|+1}(E_\omega(\z(|V_j|+1\y)^2\big|Z_0,Z_{-1},\cdots,Z_{-m+1}))^{\frac{1}{2}}\no\\
  &&\quad\quad=(|Z_{-m+1}|+1)\Big(1+2\sum_{l=1}^Lb_{-m+1}(l)\no\\
  &&\quad\quad\quad\quad\quad\quad\quad\quad\quad\quad+\sum_{l=1}^L(b_{-m+1}(l)+2b_{-m+1}^2(l))+4\sum_{1\le k<l\le
  L}b_{-m+1}(k)b_{-m+1}(l)\Big)^{\frac{1}{2}}\no\\
  &&\quad\quad=:(|Z_{-m+1}|+1)\tilde{R}( M_{-m+1})^{\frac{1}{2}}.
\end{eqnarray}
It follows that
\begin{eqnarray}
  &&E\bigg(\Big(\frac{|Z_{-m}|}{|Z_{-m+1}|+1}\Big)^\ka1_{[m\le\nu]}\bigg)
  =E\Big(E_\omega\Big(\Big(\frac{|Z_{-m}|}{|Z_{-m+1}|+1}\Big)^\ka\Big|Z_0,Z_{-1}\cdots,Z_{-m+1}\Big);m\le\nu\Big)\no\\
  &&\quad\quad\mbox{by (\ref{zmm}) and independence and stationarity}\no\\
  &&\quad\quad\le
  E(R( M_{-m+1})^{\frac{\ka}{2}};m\le\nu)=\mbe(\tilde{R}( M_{0})^\frac{\ka}{2})P(\nu>m-1).\no
\end{eqnarray}
For $1\le\ka\le2,$ since condition \textbf{(C3)} implies that
$\mbe(\tilde{R}(M_0)^\frac{\ka}{2})<\infty,$ then the inequality
(\ref{zsi}) follows from (\ref{zss}) and Theorem \ref{nu}. For
$\ka<1$ we have from (\ref{zzz}) that
\begin{eqnarray}
 && E(|Z_{-\sigma}|^\ka;\sigma<\nu)\le
 (A+1)^\ka\sum_{m\ge1}E\Big(\bigg(\frac{|Z_{-m}|}{|Z_{-m+1}|+1}\bigg)^\ka1_{[m\le\nu]}\Big)\no\\
 &&\quad\quad\le
 (A+1)^\ka\sum_{m\ge1}E\z(1_{[m\le\nu]}(|Z_{-m+1}|+1)^{-\ka}(E_\omega(|Z_{-m}|\big|Z_0,Z_{-1},\cdots,Z_{-m+1}))^\ka\y)\no\\
 &&\quad\quad\le
 (A+1)^\ka\sum_{m\ge1}E\Big(1_{[m\le\nu]}\z(1+\sum_{l=1}^Lb_{-m+1}(l)\y)^\ka\Big)\no\\
 &&\quad\quad=(A+1)^\ka\sum_{m\ge1}P(\nu>m-1)\mbe\Big(\Big(1+\sum_{l=1}^Lb_{0}(l)\Big)^\ka\Big)<\infty.\no
\end{eqnarray}
For $\ka>2,$ note that $$W=\sum_{0\le m<\nu}Y_{-m}=\sum_{m=0}^\infty
Y_{-m}1_{[m<\nu]}.$$ Then by the independence of $Y_{-m}$ and
$1_{[m<\nu]}$ we have
$$(E|W|^2)^{\frac{1}{2}}\le\sum_{m=0}^\infty(E\z(|Y_{-m}|^21_{[m<\nu]}\y))^{\frac{1}{2}}=\sum_{m=0}^\infty (E|Y_0|^2)^\frac{1}{2}(P(\nu>m))^\frac{1}{2}.
$$
Theorem \ref{nu} implies that $E|W|^2$ is finite if we can show that
$$E|Y_0|^2<\infty.$$
In fact,
\begin{eqnarray}\label{dis}
(E\left(|Y_0|^2\right))^\frac{1}{2}&=&\Big(E\bigg(\Big(\sum_{t=1}^\infty|Z(0,-t)|\Big)^2\bigg)\Big)^\frac{1}{2}\no\\
&\le&
\sum_{t=1}^\infty\z(E|Z(0,-t)|^2\y)^{\frac{1}{2}}=\sum_{t=1}^\infty\z(\mbe\left[E_\omega|Z(0,-t)|^2\right]\y)^{\frac{1}{2}}\no\\
&=&
\sum_{t=1}^\infty\Big(\mbe\left([E_\omega|Z(0,-t)|]^2\right)+\mbe\left[V_\omega(|Z(0,-t)|)\right]\Big)^{\frac{1}{2}}\no\\
&\le&
\sum_{t=1}^\infty\Big(\mbe\left([E_\omega|Z(0,-t)|]^2\right)\Big)^\frac{1}{2}+\sum_{t=1}^\infty\Big(\mbe\left[V_\omega(|Z(0,-t)|)\right]\Big)^{\frac{1}{2}}.
\end{eqnarray}To estimate the first term in the rightmost-hand side of (\ref{dis}), note that
\begin{equation}
  \mbe[(E_\omega|Z(0,-t)|)^2]=\mbe(|e_1 M_{0}M_{-1}\cdots  M_{-t+1}|^2)\le
 L^2\mbe( \parallel M_{0}M_{-1}\cdots  M_{-t+1}\parallel^2).\no
\end{equation}
Since $\ka>2,$ (\ref{wk1}) and (\ref{wk}) of Theorem \ref{kes} imply
that for some $\beta<0$
\begin{equation}\label{asy1}
\mbe( \parallel M_{0}M_{-1}\cdots  M_{-t+1}\parallel^2)\sim e^{\beta
t} \mbox{ as } t\rto.
\end{equation}
Now we estimate the second term in the rightmost-hand side of
(\ref{dis}). Note that
\begin{eqnarray}
 && V_\omega(|Z(0,-t)|)=E_\omega(\z[|Z(0,-t)|-|E_\omega(Z(0,-t))|\y]^2)\le
 E_\omega(|Z(0,-t)-E_\omega(Z(0,-t))|^2)\no\\
 &&\quad\quad=E_\omega\z(E_\omega\z[|Z(0,-t)-Z(0,-t+1) M_{-t+1}|^2\big|Z(0,-t+1)\y]\y)\no\\
  &&\quad\quad=E_\omega\z(|Z(0,-t+1)|E_\omega\z[|Z(0,-t)-Z(0,-t+1) M_{-t+1}|^2\big|Z(0,-t+1)=e_1\y]\y)\no\\
   &&\quad\quad=|e_1 M_{0}M_{-1}\cdots  M_{-t+2}|\Big(\sum_{l=1}^{L}(b_{-t+1}(l)+b_{-t+1}(l)^2)+2\sum_{1\le k<l\le
  L}b_{-t+1}(l)b_{-t+1}(k)\Big)\no
\end{eqnarray}
Then we have that
\begin{eqnarray}\label{asy2}
&&\mbe(V_\omega(|Z(0,-t)|))\no\\
&&\quad\quad\le L\mbe(\parallel M_{0}M_{-1}\cdots
M_{-t+2}\parallel)\mbe\Big(\sum_{l=1}^{L}(b_{0}(l)+b_{0}(l)^2)+2\sum_{1\le
k<l\le
  L}b_{0}(l)b_{0}(k)\Big)\no\\
  &&\quad\quad
 \le C\mbe(\parallel M_{0}M_{-1}\cdots  M_{-t+2}\parallel)\sim
 Ce^{\gamma
 t}
\end{eqnarray}
 for some $C>0$ and $\gamma<0$ for the same reason as (\ref{asy1}).
 Then (\ref{asy1}) together with  (\ref{asy2}) implies the
 convergence of the series in the rightmost-hand side of (\ref{dis}). Therefore
 $E|Y_0|^2<\infty.$
 \qed

Recall  $$S_{-\sigma,-m}=\mbox{number of progeny alive at time $-m$
of the $Z_{-\sigma}$ particles present at $-\sigma>-m$.}$$ Let
$S_{-\sigma,-\sigma}=Z_{-\sigma},$  and
$$S_{-\sigma}=\sum_{m=\sigma}^\infty S_{-\sigma,-m}=Z_{-\sigma}+\mbox{total progeny of the $Z_{-\sigma}$ particles at $-\sigma.$}$$

\noindent{\bf Proof of Lemma \ref{sz}:} Observe that
\[S_{-\sigma,-m}-Z_{-\sigma}\prod_{i=\sigma}^{m-1} M_{-i}=\sum_{\sigma+1\le l\le
m}\Big(S_{-\sigma,-l}\prod_{i=l}^{m-1}
M_{-i}-S_{-\sigma,-l+1}\prod_{i=l-1}^{m-1} M_{-i}\Big),\] the
convention being that empty product is $I,$ and therefore
\begin{eqnarray}
  &&\sum_{m=\sigma}^\infty\Big(S_{-\sigma,-m}-Z_{-\sigma}\prod_{i=\sigma}^{m-1} M_{-i}\Big)=\sum_{l=\sigma+1}^\infty\sum_{m=l}^\infty\Big(S_{-\sigma,-l}\prod_{i=l}^{m-1} M_{-i}-S_{-\sigma,-l+1}\prod_{i=l-1}^{m-1} M_{-i}\Big)\no\\
  &&\quad\quad\quad\quad\quad\quad=\sum_{l=\sigma+1}^\infty(S_{-\sigma,-l}-S_{-\sigma,-l+1} M_{-l+1})\sum_{m=l}^\infty\prod_{i=l}^{m-1} M_{-i}\no\\
  &&\quad\quad\quad\quad\quad\quad=\sum_{l=\sigma+1}^\infty(S_{-\sigma,-l}-S_{-\sigma,-l+1} M_{-l+1})(I+\eta_{-l}).\no
\end{eqnarray}
Then we have
\begin{eqnarray}
  &&P\Big(\Big|\sum_{m=\sigma}^\infty\Big(S_{-\sigma,-m}-Z_{-\sigma}\prod_{i=\sigma}^{m-1} M_{-i}\Big)\Big|\ge \epsilon
  x,\sigma<\nu\Big)\no\\
   &&\hspace{1cm}=\sum_{j=1}^\infty P\Big(\Big|\sum_{m=\sigma}^\infty\Big(S_{-\sigma,-m}-Z_{-\sigma}\prod_{i=\sigma}^{m-1} M_{-i}\Big)\Big|\ge \epsilon
  x,\sigma<\nu,\sigma=j\Big)\no\\
  &&\hspace{1cm}=\sum_{j=1}^\infty P\Big(\Big|\sum_{l=\sigma+1}^\infty(S_{-\sigma,-l}-S_{-\sigma,-l+1} M_{-l+1})(I+\eta_{-l})\Big|\ge \epsilon
  x,\sigma<\nu,\sigma=j\Big)\no\\
&&\hspace{1cm}\le \sum_{j=1}^\infty
P\Big(\sum_{l=\sigma+1}^\infty\Big|(S_{-\sigma,-l}-S_{-\sigma,-l+1}
M_{-l+1})\Big||\textbf{1}(I+\eta_{-l})|\ge \epsilon
  x,\sigma<\nu,\sigma=j\Big).\no
\end{eqnarray}
By a similar argument as (\ref{duli}), the rightmost-hand side of
the last expression is less than or equal to
\begin{eqnarray}
&&\sum_{j=1}^\infty\sum_{l=j+1}^\infty
P(|S_{-\sigma,-l}-S_{-\sigma,-l+1} M_{-l+1}|e_0(I+\eta_{-l})e_0^T\ge
\frac{1}{2}(l-\sigma)^{-2}L^{-2}\epsilon
  x,\sigma<\nu,\sigma=j)\no\\
  &&\hspace{1cm}=\sum_{j=1}^\infty\sum_{l=j+1}^\infty\int P(|S_{-j,-l}-S_{-j,-l+1} M_{-l+1}|\in ds,j<\nu,\sigma=j,\no\\
  &&\hspace{5cm} e_0(I+\eta_{-l})e_0^T\ge
(2s)^{-1}(l-j)^{-2}L^{-2}\epsilon
  x).\no
\end{eqnarray}
Note that $\{|S_{-j,-l}-S_{-j,-l+1} M_{-l+1}|\in
ds,j<\nu,\sigma=j\},$ defined in term of
$\{Z_0,Z_{-1},...,Z_{-l}\},$ depends only on
$\sigma(\omega_{-i};i<l)$ and that $\{e_0(I+\eta_{-l})e_0^T\ge
(2s)^{-1}(l-j)^{-2}L^{-2}\epsilon
  x\}$ depends only on $\sigma(\omega_{-i};i\ge l).$
Then by the independence and stationarity of the environment, the
right-hand side of the above equality equals to
\begin{eqnarray}
&&\sum_{j=1}^\infty\sum_{l=j+1}^\infty\int
P(|S_{-j,-l}-S_{-j,-l+1} M_{-l+1}|\in ds,j<\nu,\sigma=j)\no\\
  &&\hspace{4cm}\times P(e_0(I+\eta_0)e_0^T\ge
(2s)^{-1}(l-j)^{-2}L^{-2}\epsilon
  x)\no
\end{eqnarray}

 which, by Theorem \ref{kes}  with $K_7$ as
in (\ref{k7}), is less than or equal to
\begin{eqnarray}
  &&x^{-\ka}(\frac{2}{\epsilon})^\ka
  K_7L^{2\ka}\sum_{j=1}^\infty\sum_{l=j+1}^\infty (l-j)^{2\ka}\int s^\ka P(|S_{-j,-l}-S_{-j,-l+1} M_{-l+1}|\in
  ds,j<\nu,\sigma=j)\no\\
  &&\quad\quad= x^{-\ka}(\frac{2}{\epsilon})^\ka
  K_7L^{2\ka}\sum_{j=1}^\infty\sum_{l=j+1}^\infty(l-j)^{2\ka}
  E(|S_{-j,-l}-S_{-j,-l+1} M_{-l+1}|^\ka;j<\nu,\sigma=j)\no\\
&&\quad\quad= x^{-\ka}(\frac{2}{\epsilon})^\ka
  K_7L^{2\ka}
  E\Big(\sum_{l=\sigma+1}^\infty(l-\sigma)^{2\ka}|S_{-\sigma,-l}-S_{-\sigma,-l+1} M_{-l+1}|^\ka;\sigma<\nu\Big)\no\\
  &&\quad\quad= x^{-\ka}(\frac{2}{\epsilon})^\ka
  K_7L^{2\ka}
  E\Big(\sum_{l=\sigma+1}^\infty(l-\sigma)^{2\ka}\no\\
  &&\hspace{3cm}\times E\Big(|S_{-\sigma,-l}-S_{-\sigma,-l+1}
   M_{-l+1}|^\ka;\sigma<\nu\Big|\omega,
  \sigma,Z_0,Z_{-1},\cdots,Z_{-\sigma}\Big)\Big).\no
\end{eqnarray}
Recalling that $\ka\le 2,$  Jensen's inequality implies that the
rightmost-hand side of above expression is less than or equal to
\begin{eqnarray}\label{eps}
&&x^{-\ka}(\frac{2}{\epsilon})^\ka
  K_7L^{2\ka}
  E\Big(\sum_{l=\sigma+1}^\infty(l-\sigma)^{2\ka}\times\no\\
  &&\hspace{.5cm}\Big\{E\Big(|S_{-\sigma,-l}-S_{-\sigma,-l+1}
   M_{-l+1}|^2;\sigma<\nu\Big|\omega,
  \sigma,Z_0,Z_{-1},\cdots,Z_{-\sigma}\Big)\Big\}^{\frac{\ka}{2}}\Big).
\end{eqnarray}

 Again as in (\ref{fc1}) we have
\begin{eqnarray}
&&E\z(|S_{-\sigma,-l}-S_{-\sigma,-l+1} M_{-l+1}|^2\big|\omega,\sigma,Z_0,Z_{-1},\cdots,Z_{-\sigma}, S_{-\sigma,-l+1}\y)\no\\
&&\quad\quad=|S_{-\sigma,-l+1}|\Big(\sum_{j=1}^L(b_{-l+1}(j)+b_{-l+1}^2(j))+2\sum_{1\le
i<j\le L}b_{-l+1}(i)b_{-l+1}(j)\Big)\no\\
&&\quad\quad=:|S_{-\sigma,-l+1}|R( M_{-l+1})\no
\end{eqnarray}
and
\begin{eqnarray}
&&(E\z(|S_{-\sigma,-l}-S_{-\sigma,-l+1} M_{-l+1}|^2\big|\omega,\sigma,Z_0,Z_{-1},\cdots,Z_{-\sigma} \y))^{\frac{\ka}{2}}\no\\
&&\quad\quad=(E\z(|S_{-\sigma,-l+1}|
\big|\omega,\sigma,Z_0,Z_{-1},\cdots,Z_{-\sigma}\y))^{\frac{\ka}{2}}
R( M_{-l+1})^{\frac{\ka}{2}}\no\\
&&\quad\quad=\Big|Z_{-\sigma}\prod_{i=\sigma}^{l-2}
M_{-i}\Big|^{\frac{\ka}{2}}R( M_{-l+1})^\frac{\ka}{2}.\no
\end{eqnarray}
Substituting to (\ref{eps}), we get that
\begin{eqnarray}
&&P\Big(\Big|\sum_{m=\sigma}^\infty\Big(S_{-\sigma,-m}-Z_{-\sigma}
\prod_{i=\sigma}^{m-1} M_{-i}\Big)\Big|\ge
\epsilon x;\sigma<\nu\Big)\no\\
 &&\quad\quad\le x^{-\ka}(\frac{2}{\epsilon})^\ka
K_7L^{2\ka}E\Big(\sum_{l=\sigma+1}^\infty(l-\sigma)^{2\ka}\Big|Z_{-\sigma}\prod_{i=\sigma}^{l-2} M_{-i}\Big|^{\frac{\ka}{2}}R( M_{-l+1})^\frac{\ka}{2};\sigma<\nu\Big)\no\\
&&\quad\quad\le x^{-\ka}(\frac{2}{\epsilon})^\ka
K_7L^{2\ka}\sum_{m=1}^\infty\sum_{l=m+1}^\infty
(l-m)^{2\ka}E\Big(|Z_{-m}|^\frac{\ka}{2}\Big\|\prod_{i=m}^{l-2}
M_{-i}\Big\|^{\frac{\ka}{2}}R(
M_{-l+1})^{\frac{\ka}{2}};m=\sigma<\nu\Big).\no
\end{eqnarray}
Again, using independence and stationarity, the rightmost-hand side
of the above expression
\begin{eqnarray}
&&\quad\quad= x^{-\ka}(\frac{2}{\epsilon})^\ka
K_7L^{2\ka}\no\\
&&\quad\quad\quad\times\sum_{m=1}^\infty\sum_{l=m+1}^\infty
(l-m)^{2\ka}\mbe\Big(\Big\|\prod_{i=0}^{l-m-2}
M_{-i}\Big\|^{\frac{\ka}{2}}\Big)\mbe\z(R( M_{0})^{\frac{\ka}{2}}\y)
 E\Big(|Z_{-m}|^\frac{\ka}{2};m=\sigma<\nu\Big)\no\\
 &&\quad\quad= x^{-\ka}(\frac{2}{\epsilon})^\ka
K_7L^{2\ka}\no\\
&&\quad\quad\quad\times\sum_{m=1}^\infty\sum_{l=1}^\infty
l^{2\ka}\mbe\Big(\Big\|\prod_{i=0}^{l-2}
M_{-i}\Big\|^{\frac{\ka}{2}}\Big)\mbe\z(R( M_{0})^{\frac{\ka}{2}}\y)
 E\Big(|Z_{-m}|^\frac{\ka}{2};m=\sigma<\nu\Big)\no\\
 && \hspace{2cm}\mbox{using 1) and 2) of Theorem \ref{kes} and condition
 \textbf{(C3)}}\no\\
&&\quad\quad\le K_8(\epsilon
x)^{-\ka}E(|Z_{-\sigma}|^{\frac{\ka}{2}};\sigma<\nu)\le K_8(\epsilon
x)^{-\ka}A^{-\frac{\ka}{2}}E(|Z_{-\sigma}|^\ka;\sigma<\nu)\no\\
&&\quad\quad\le \epsilon x^{-\ka}E(|Z_{-\sigma}|^\ka;\sigma<\nu)\no
\end{eqnarray}
for $A\ge A_1(\epsilon),$ and some $K_8>0.$ \qed

\noindent{\large{\bf Acknowledgements:}} The authors would like to
thank Professor Zenghu Li for his useful discussions on the
multitype branching processes, we also thank Ms. Hongyan Sun and Lin
Zhang for the stimulating discussions. The project is partially
supported by National Nature Science Foundation of China (Grant
No.10721091) and NCET( No.\,05-0143). 

{\center

\end{document}